\documentclass{amsart}
\usepackage{amssymb, latexsym}
\usepackage{enumerate}

\newcommand\dive{\operatorname{div}}

\newcommand\card{\operatorname{card}}
\newcommand\gu{\operatorname{g(u)}}

\newcommand\gfivesixM{\operatorname{g^{\frac{5}{6}}(M) }}
\newcommand\gfivesixu{\operatorname{g^{\frac{5}{6}} (u)}}
\newcommand\gonesixM{\operatorname{g^{\frac{1}{6}} (M) }}
\newcommand\gonesixu{\operatorname{g^{\frac{1}{6}} (u) }}
\newcommand\gfivetwfourM{\operatorname{g^{\frac{5}{24}}(M) }}

\theoremstyle{plain}
\newtheorem{thm}{Theorem}

\newtheorem{rem}[thm]{Remark}
\newtheorem{prop}{Proposition}
\newtheorem{lem}[prop]{Lemma}

\newtheorem{res}{Result}

\numberwithin{equation}{section} \numberwithin{thm}{section}

\begin{document}

\title[Loglog Supercritical Wave Equation ] {Global Existence Of Smooth Solutions Of A $3D$ Loglog
energy-supercritical Wave Equation}
\author{Tristan Roy}
\address{University of California, Los Angeles}
\email{triroy@math.ucla.edu}

\begin{abstract}
We prove global existence of smooth solutions of the $3D$ loglog
energy-supercritical wave equation $\partial_{tt} u - \triangle u =
-u^{5} \log^{c} \left( log(10+u^{2}) \right) $ with $0 < c <
\frac{8}{225}$ and smooth initial data $(u(0)=u_{0}, \, \partial_{t}
u(0)=u_{1})$. First we control the $L_{t}^{4} L_{x}^{12}$ norm of
the solution on an arbitrary size time interval by an expression
depending on the energy and an \textit{a priori} upper bound of its
$L_{t}^{\infty} \tilde{H}^{2}(\mathbb{R}^{3})$ norm, with
$\tilde{H}^{2}(\mathbb{R}^{3}):=\dot{H}^{2}(\mathbb{R}^{3}) \cap
\dot{H}^{1}(\mathbb{R}^{3})$. The proof of this long time estimate
relies upon the use of some potential decay estimates \cite{bahger,
shatstruwe} and a modification of an argument in \cite{taolog}. Then
we find an \textit{a posteriori} upper bound of the $L_{t}^{\infty}
\tilde{H}^{2}(\mathbb{R}^{3})$ norm of the solution by combining the
long time estimate with an induction on time of the Strichartz
estimates.
\end{abstract}

\maketitle

\section{Introduction}

We shall consider the defocusing loglog energy-supercritical wave equation

\begin{equation}
\begin{array}{ll}
\partial_{tt} u - \triangle u & = -f(u)
\end{array}
\label{Eqn:LogQuint}
\end{equation}
where $u: \mathbb{R} \times \mathbb{R}^{3} \rightarrow \mathbb{R}$
is a real-valued scalar field and $f(u):= u^{5} g(u)$ with
$g(u):=\log^{c} \left( log(10+u^{2}) \right)$, $ 0 < c <
\frac{8}{225}$. \textit{Classical solutions} of (\ref{Eqn:LogQuint})
are solutions that are infinitely differentiable and compactly
supported in space for each fixed time $t$. It is not difficult to
see that \textit{classical solutions} of (\ref{Eqn:LogQuint})
satisfy the energy conservation law

\begin{equation}
\begin{array}{ll}
E & := \frac{1}{2} \int_{\mathbb{R}^{3}}  \left( \partial_{t} u (t,x) \right)^{2} \, dx + \frac{1}{2} \int_{\mathbb{R}^{3}} |\nabla u (t,x)|^{2}
\, dx + \int_{\mathbb{R}^{3}} F(u(t,x)) \, dx
\end{array}
\label{Eqn:DefE}
\end{equation}
where $F(u):= \int_{0}^{u} f(v) \, dv $. \textit{Classical solutions} of (\ref{Eqn:LogQuint}) enjoy three symmetry properties that we use
throughout this paper

\begin{itemize}

\item the \textit{ time translation invariance}: if $u$ is a solution of (\ref{Eqn:LogQuint}) and $t_{0}$ is a fixed time then
$\tilde{u}(t,x):=u(t-t_{0},x)$ is also a solution of
(\ref{Eqn:LogQuint}).

\item the \textit{space translation invariance}: if $u$ is a solution of (\ref{Eqn:LogQuint}) and $x_{0}$ is a fixed point lying in
$\mathbb{R}^{3}$ then $\tilde{u}(t,x):=u(t,x-x_{0})$ is also a
solution of (\ref{Eqn:LogQuint})

\item the \textit{time reversal invariance}: if $u$ is a solution to (\ref{Eqn:LogQuint}) then $\tilde{u}(t,x):=u(-t,x)$ is also a solution to
(\ref{Eqn:LogQuint})
\end{itemize}

The defocusing loglog energy-supercritical wave equation (\ref{Eqn:LogQuint}) is closely related to the power-type defocusing wave equations,
namely

\begin{equation}
\begin{array}{ll}
\partial_{tt} u - \triangle u & = - |u|^{p-1} u
\end{array}
\label{Eqn:PDef}
\end{equation}
Solutions of (\ref{Eqn:PDef}) have an invariant scaling

\begin{equation}
\begin{array}{ll}
u(t,x) & \rightarrow u^{\lambda}(t,x):=\frac{1}{\lambda^{\frac{2}{p-1}}} u \left( \frac{t}{\lambda}, \, \frac{x}{\lambda} \right)
\end{array}
\end{equation}
and (\ref{Eqn:PDef}) is $s_{c}:= \frac{3}{2} - \frac{2}{p-1}$
critical. This means that the $\dot{H}^{s_{c}}(\mathbb{R}^{3})
\times \dot{H}^{s_{c}-1}(\mathbb{R}^{3})$ norm of $\left( u(0),
\partial_{t} u(0) \right)$ is invariant under scaling or, in other
words, $\| u^{\lambda}(0) \|_{\dot{H}^{s_{c}}(\mathbb{R}^{3})} = \|
u(0) \|_{\dot{H}^{s_{c}}(\mathbb{R}^{3})}$ and $\|
\partial_{t} u^{\lambda}(0) \|_{\dot{H}^{s_{c}-1}(\mathbb{R}^{3})} = \| \partial_{t}
u (0) \|_{\dot{H}^{s_{c}-1}(\mathbb{R}^{3})}$. If $p=5$ then
$s_{c}=1$ and this is why the quintic defocusing cubic wave equation

\begin{equation}
\begin{array}{ll}
\partial_{tt} u -\triangle u & = - u^{5}
\end{array}
\label{Eqn:QuintDefocWave}
\end{equation}
is called the energy-critical equation. If $1<p<5$ then $s_{c} < 1$ and (\ref{Eqn:PDef}) is energy-subcritical while if $p>5$ then $s_{c}> 1$
and (\ref{Eqn:PDef}) is energy-supercritical. Notice that for every $p>5$ there exists two positive constant $\lambda_{1}(p)$, $\lambda_{2}(p)$
such that

\begin{equation}
\begin{array}{lll}
\lambda_{1}(p)|u|^{5} \leq & |f(u)| & \leq \lambda_{2}(p) \max{ \left( 1,|u|^{p} \right)}
\end{array}
\end{equation}
This is why (\ref{Eqn:LogQuint}) is said to belong to the group of barely supercritical equations. There is another way to see that. Notice that
a simple integration by part shows that

\begin{equation}
\begin{array}{ll}
F(u) & \sim \frac{u^{6}}{6} g(u)
\end{array}
\label{Eqn:AsympF}
\end{equation}
and consequently the nonlinear potential term of the energy $\int_{\mathbb{R}^{3}} F(u) \, dx \sim \int_{\mathbb{R}^{3}} u^{6} g(u) \, dx$ just
barely fails to be controlled by the linear component, in contrast to (\ref{Eqn:QuintDefocWave}).

The energy-critical wave equation (\ref{Eqn:QuintDefocWave}) has received a great deal of attention. Grillakis \cite{grill,grill2} established
global existence of smooth solutions (global regularity) of this equation with smooth initial data $(u(0)=u_{0},\partial_{t} u(0)=u_{1})$. His
work followed the work of Rauch for small data \cite{rauch} and the one of Struwe \cite{St} handling the spherically symmetric case. Later
Shatah and Struwe \cite{shatstruwe} gave a simplified proof of this result. Kapitanski \cite{kap} and, independently, Shatah and Struwe
\cite{shatstruwe2}, proved global existence of solutions with data $(u_{0},u_{1})$ in the energy class.

We are interested in proving global regularity of
(\ref{Eqn:LogQuint}) with smooth initial data $(u_{0}, u_{1}) $. By
standard persistence of regularity results it suffices to prove
global existence of solutions $u \in  \mathcal{C} \left( [0, \, T],
\, \tilde{H}^{2}(\mathbb{R}^{3}) \right) \cap \mathcal{C}^{1} \left(
[0, \, T], H^{1}(\mathbb{R}^{3})  \right) $ with data $(u_{0},
u_{1}) \in \tilde{H}^{2} (\mathbb{R}^{3}) \times
H^{1}(\mathbb{R}^{3})$. Here $\tilde{H}^{2}(\mathbb{R}^{3})$ denotes
the following space

\begin{equation}
\begin{array}{ll}
\tilde{H}^{2}(\mathbb{R}^{3}) & : = \dot{H}^{2}(\mathbb{R}^{3}) \cap
\dot{H}^{1}(\mathbb{R}^{3})
\end{array}
\end{equation}
In view of the local well-posedness theory \cite{linsog}, standard
limit arguments and the finite speed of propagation it suffices to
find an \textit{a priori} upper bound of the form

\begin{equation}
\begin{array}{ll}
\left\| ( u(T), \partial_{t} u (T) ) \right\|_{\tilde{H}^{2}
(\mathbb{R}^{3}) \times H^{1}(\mathbb{R}^{3})} & \leq C_{1} \left(
\| u_{0} \|_{\tilde{H}^{2}(\mathbb{R}^{3})}, \, \| u_{1}
\|_{H^{1}(\mathbb{R}^{3})}, \, T   \right)
\end{array}
\end{equation}
for all times $T>0$ and for \textit{classical solutions} $u$ of
(\ref{Eqn:LogQuint}) with smooth and compactly supported data
$(u_{0},u_{1})$. Here $C_{1}$ is a constant depending only on $\|
u_{0} \|_{\tilde{H}^{2}(\mathbb{R}^{3})}$, $\| u_{1}
\|_{H^{1}(\mathbb{R}^{3})}$ and the time $T$.

The global behaviour of the solutions of the supercritical wave
equations is poorly understood, mostly because of the lack of
conservation laws in $\tilde{H}^{2}(\mathbb{R}^{3})$. Nevertheless
Tao \cite{taolog} was able to prove global regularity for another
barely supercritical equation, namely

\begin{equation}
\begin{array}{ll}
\partial_{tt} u - \triangle u & = - u^{5} \log{(2+u^{2})}
\end{array}
\end{equation}
with radial data. The main result of this paper is

\begin{thm}
The solution of (\ref{Eqn:LogQuint}) with smooth data $(u_{0},
u_{1})$ exists for all time. Moreover there exists a nonnegative
constant $ M_{0}=M_{0} ( \| u_{0}
\|_{\tilde{H}^{2}(\mathbb{R}^{3})}, \| u_{1}
\|_{H^{1}(\mathbb{R}^{3} )})$ depending only on $\| u_{0}
\|_{\tilde{H}^{2}(\mathbb{R}^{3})}$ and $\| u_{1}
\|_{H^{1}(\mathbb{R}^{3})}$ such that

\begin{equation}
\begin{array}{ll}
\| u \|_{L_{t}^{\infty} \tilde{H}^{2} (\mathbb{R} \times
\mathbb{R}^{3})} + \| \partial_{t} u  \|_{L_{t}^{\infty} H^{1}
(\mathbb{R} \times \mathbb{R}^{3})} & \leq M_{0}
\end{array}
\end{equation}
\label{Thm:GlobLogQuint}
\end{thm}

We recall some basic properties and estimates. Let $Q$ be a function, let $J$ be an interval and let $t_{0} \in J$ be a fixed time. If $u$ is a
\textit{classical solution} of the more general problem $\partial_{tt} u - \triangle u = Q$ then $u$ satisfies the Duhamel formula

\begin{equation}
\begin{array}{ll}
u(t) & = u_{l,t_{0}}(t) + u_{nl,t_{0}}(t), \, t \in J
\end{array}
\end{equation}
with $u_{l,t_{0}}$, $u_{nl,t_{0}}$ denoting the linear part and the nonlinear part respectively of the solution starting from $t_{0}$. Recall
that

\begin{equation}
\begin{array}{ll}
u_{l,t_{0}}(t) & = \cos{(t-t_{0})D} u(t_{0}) + \frac{\sin{(t-t_{0})D}}{D} \partial_{t} u(t_{0})
\end{array}
\end{equation}
and

\begin{equation}
\begin{array}{ll}
u_{nl,t_{0}}(t) &  = - \int_{t_{0}}^{t} \frac{\sin{(t-t^{'})D}}{D}
Q(t^{'}) \, dt^{'}
\end{array}
\end{equation}
with $D$ the multiplier defined by $\widehat{Df}(\xi):= |\xi|
\widehat{f}(\xi)$. An explicit formula for $
\frac{\sin{(t-t^{'})D}}{D} Q(t^{'})$ and $t \neq t^{'}$ is

\begin{equation}
\begin{array}{ll}
 \left[ \frac{\sin{(t-t^{'}) D}}{D} Q(t^{'}) \right] (x) &  = \frac{1}{4 \pi |t- t^{'}|} \int_{|x- x^{'}|=|t-t^{'}|} Q(t^{'},x^{'}) \, dS(x^{'})
\end{array}
\label{Eqn:ExplicitCompSin}
\end{equation}
For a proof see for example \cite{sogge}. We recall that $u_{l,t_{0}}$  satisfies $\partial_{tt} u_{l,t_{0}} - \triangle u_{l,t_{0}}  = 0$, $
u_{l,t_{0}}(t_{0})= u(t_{0})$ and $ \partial_{t} u_{l,t_{0}}(t_{0}) = \partial_{t} u(t_{0}) $ while $u_{nl,t_{0}}$ is the solution of $
\partial_{tt} u_{nl,t_{0}} - \triangle  u_{nl,t_{0}} = Q $, $  u_{nl,t_{0}}(t_{0}) = 0 $ and  $ \partial_{t}  u_{nl,t_{0}}(t_{0}) = 0 $. We
recall the Strichartz estimate \cite{ginebvelo,keeltao,linsog,sogge}

\begin{equation}
\begin{array}{ll}
\| u \|_{L_{t}^{q} L_{x}^{r} (J \times \mathbb{R}^{3})}  & \lesssim \| \partial_{t} u (t_{0}) \|_{L_{x}^{2}(\mathbb{R}^{3})} + \| \nabla
u(t_{0}) \|_{L_{x}^{2} (\mathbb{R}^{3})}   + \| Q \|_{L_{t}^{1} L_{x}^{2} (J \times \mathbb{R}^{3}) }
\end{array}
\label{Eqn:Str}
\end{equation}
if $(q,r)$ is wave admissible, i.e $(q,r) \in (2, \, \infty] \times
[2, \, \infty]$ and $\frac{1}{q} + \frac{3}{r}=\frac{1}{2}$.

We set some notation that appear throughout the paper. We write
$C=C(a_{1},...,a_{n})$ if $C$ only depends on the parameters
$a_{1}$,...,$a_{n}$. We write $A \lesssim B$ if there exists a
universal nonnegative constant $C^{'}>0$ such that $A \leq C^{'} B$.
$A=O(B)$ means $A \lesssim B$. More generally we write $A
\lesssim_{a_{1},....,a_{n}} B$ if there exists a nonnegative
constant $C^{'}=C(a_{1},...,a_{n})$ such that $A \leq C^{'} B$. We
say that $C^{''}$ is the constant determined by $\lesssim$ in $A
\lesssim_{a_{1},...,a_{n}} B $ if $C^{''}$ is the smallest constant
among the $C^{'}$ s such that $A \leq C^{'}B$. We write $A
<<_{a_{1},..,a_{n}} B$ if there exists a universal nonnegative small
constant $c=c(a_{1},...,a_{n})$ such that $A \leq c B$. Similar
notions are defined for $A \gtrsim  B$, $A \gtrsim_{a_{1},...,a_{n}}
B$ and $A >> B$. In particular we say that $C^{''}$ is the constant
determined by $\gtrsim$ in $A \gtrsim B$ if $C^{''}$ is the largest
constant among the $C^{'}$s such that $A \geq C^{'} B$. If $x$ is
number then $x+$ and $x-$ are slight variations of $x$: $x+ := x +
\alpha \epsilon$ and  $x- := x - \beta \epsilon $ for some $\alpha
> 0$, $\beta >0$ and $0< \epsilon << 1$.

Let $\Gamma_{+}$ denote the forward light cone

\begin{equation}
\begin{array}{ll}
\Gamma_{+} & = \left\{ (t,x): t>|x| \right\}
\end{array}
\end{equation}
and if $J=[a,b]$ is an interval let $\Gamma_{+}(J)$ denote the light cone truncated to $J$ i.e

\begin{equation}
\begin{array}{ll}
\Gamma_{+}(J) & := \Gamma_{+} \cap (J \times \mathbb{R}^{3})
\end{array}
\end{equation}
Let $e(t)$ denote the local energy i.e

\begin{equation}
\begin{array}{ll}
e(t) & : = \frac{1}{2} \int_{|x| \leq t} \left( \partial_{t} u (t,x) \right)^{2} \, dx + \frac{1}{2} \int_{|x| \leq t} \left| \nabla u(t,x)
\right|^{2} \, dx + \int_{|x| \leq t} F ( u(t,x) ) \, dx
\end{array}
\end{equation}
If $u$ is a solution of (\ref{Eqn:LogQuint})  then by using the finite speed of propagation and the Strichartz estimates we have

\begin{equation}
\begin{array}{ll}
\| u \|_{L_{t}^{q} L_{x}^{r}(\Gamma_{+}(J))} & \lesssim \| \nabla u(b) \|_{L_{x}^{2}(\mathbb{R}^{3})} +
\| \partial_{t} u(b) \|_{L_{x}^{2}(\mathbb{R}^{3})} + \| Q \|_{L_{t}^{1} L_{x}^{2} (\Gamma_{+}(J))}
\end{array}
\label{Eqn:StrCone}
\end{equation}
if $(q,r)$ is wave admissible. If $J_{1}:=[a_{1},a_{2}]$ and $J_{2}:=[a_{2},a_{3}]$ then we also have

\begin{equation}
\begin{array}{ll}
\| u \|_{L_{t}^{q} L_{x}^{r} (\Gamma_{+} (J_{1}))} & \lesssim  \| \nabla u(a_{3}) \|_{L_{x}^{2}(\mathbb{R}^{3})} + \| \partial_{t} u(a_{3})
\|_{L_{x}^{2}(\mathbb{R}^{3})} + \| Q \|_{L_{t}^{1} L_{x}^{2} (\Gamma_{+}(J_{1} \cup J_{2}))}
\end{array}
\label{Eqn:Strcone2}
\end{equation}
We recall also the well-known Sobolev embeddings. If $h$ is a smooth
function then

\begin{equation}
\begin{array}{ll}
\| h \|_{L^{\infty}(\mathbb{R}^{3})} & \lesssim \| h
\|_{\tilde{H}^{2}(\mathbb{R}^{3})}
\end{array}
\label{Eqn:SobEmbed}
\end{equation}
and

\begin{equation}
\begin{array}{ll}
\| h \|_{L^{6}(\mathbb{R}^{3})} & \lesssim  \| \nabla h \|_{L^{2}(\mathbb{R}^{3})}
\end{array}
\label{Eqn:SobEmbed2}
\end{equation}
If $u$ is the solution of (\ref{Eqn:LogQuint}) with data $(u_{0},
u_{1}) \in \tilde{H}^{2}(\mathbb{R}^{3}) \times
H^{1}(\mathbb{R}^{3})$ then we get from (\ref{Eqn:SobEmbed})

\begin{equation}
\begin{array}{ll}
E & \lesssim \| u_{0} \|_{\tilde{H}^{2} (\mathbb{R}^{3})}^{2} \max
\left( 1, \| u_{0} \|^{4}_{\tilde{H}^{2}(\mathbb{R}^{3})} g( \|
u_{0} \|_{\tilde{H}^{2}(\mathbb{R}^{3})} ) \right)
\end{array}
\label{Eqn:BoundE}
\end{equation}
We shall use the Paley-Littlewood technology. Let $\phi(\xi)$ be a
bump function adapted to $\left\{ \xi \in \mathbb{R}^{3}: \, |\xi|
\leq  2 \right\}$  and equal to one on $\left\{ \xi \in
\mathbb{R}^{3}: \, |\xi| \leq 1 \right\}$. If $(M,N) \in
2^{\mathbb{Z}} \times 2^{\mathbb{Z}}$ are dyadic numbers then the
Paley-Littlewood projection operators $P_{M}$, $P_{<N}$ and $P_{\geq
N}$ are defined in the Fourier domain by

\begin{equation}
\begin{array}{ll}
\widehat{P_{M} f}(\xi) : = \left( \phi \left( \frac{\xi}{M} \right)
- \phi \left( \frac{\xi}{2M} \right) \right) \hat{f}(\xi)
\end{array}
\end{equation}

\begin{equation}
\begin{array}{ll}
\widehat{P_{<N} f}(\xi) : = \sum \limits_{M < N }{} \widehat{P_{M}
f}(\xi)
\end{array}
\end{equation}
and

\begin{equation}
\begin{array}{ll}
\widehat{P_{\geq  N} f}(\xi) : = \sum \limits_{M \geq N }{}
\widehat{P_{M} f}(\xi)
\end{array}
\end{equation}
The inverse Sobolev inequality can be stated as follows

\begin{prop}{\textbf{"Inverse Sobolev inequality"} \cite{tao}}
Let $g$ be a smooth function such that $\| g \|_{\dot{H}^{1}
(\mathbb{R}^{3})} \lesssim E^{\frac{1}{2}}$ and $\| P_{\geq N} g
\|_{L_{x}^{6}(\mathbb{R}^{3})} \gtrsim \eta$ for some real number
$\eta > 0$ and for some dyadic number $N>0$. Then there exists a
ball $B(x,r) \subset \mathbb{R}^{3}$ with $r= O \left( \frac{1}{N}
\right) $ such that we have the mass concentration estimate

\begin{equation}
\begin{array}{ll}
\int_{B(x,r)} |g(y)|^{2} \, dy & \gtrsim \eta^{3} E^{-\frac{1}{2}}
r^{2}
\end{array}
\end{equation}
\label{prop:InvSob}
\end{prop}
We also recall a result that shows that the mass of solutions of
(\ref{Eqn:LogQuint}) can be locally in time controlled

\begin{prop}{\textbf{"Local mass is locally stable"} \cite{tao}}
Let $J$ be a time interval, let $t$, $t^{'} \in J$ and let $B(x,r)$ be a ball. Let $u$ be a solution of (\ref{Eqn:LogQuint}). Then

\begin{equation}
\begin{array}{ll}
 \left(  \int_{B(x,r)} |u(t^{'},y)|^{2} \, dy \right)^{\frac{1}{2}} & = \left( \int_{B(x,r)} |u(t,y)|^{2} \, dy   \right)^{\frac{1}{2}}
 +O \left( E^{\frac{1}{2}} |t - t^{'}| \right)
\end{array}
\end{equation}
\label{prop:tao}
\end{prop}
Notice that this result, proved for (\ref{Eqn:QuintDefocWave}) in
\cite{tao}, is also true for (\ref{Eqn:LogQuint}). Indeed the proof
relied upon the fact that the $L^{2}(\mathbb{R}^{3})$ norm of the
velocity of the solution of (\ref{Eqn:QuintDefocWave}) at time $t$
is bounded by the square root of its energy, which is also true for
the solution of (\ref{Eqn:LogQuint}) ( by (\ref{Eqn:DefE}) and
(\ref{Eqn:AsympF}) ).

Now we make some comments with respect to Theorem
\ref{Thm:GlobLogQuint}. If the function $g$ were a positive
constant, then it would be easy to prove that the solution of
(\ref{Eqn:LogQuint}) with data $(u_{0},u_{1}) \in
\tilde{H}^{2}(\mathbb{R}^{3}) \times H^{1}(\mathbb{R}^{3})$, since
we have a good global theory for (\ref{Eqn:QuintDefocWave}).
Therefore we can hope to prove global well-posedness for $g$ slowly
increasing to infinity, by extending the technology to prove global
well-posedness for (\ref{Eqn:QuintDefocWave}). Notice also that Tao,
in \cite{tao}, found that the solution $u$ of
(\ref{Eqn:QuintDefocWave}) satisfies

\begin{equation}
\begin{array}{ll}
\| u \|_{L_{t}^{4} L_{x}^{12}(\mathbb{R} \times \mathbb{R}^{3})} &
\lesssim \tilde{E}^{\tilde{E}^{O(1)}}
\end{array}
\label{Eqn:BoundEEnerg}
\end{equation}
with $\tilde{E}$ the energy of $u$. The structure of $g$ is a double
log: it is, roughly speaking, the inverse function of the towel
exponential bound in (\ref{Eqn:BoundEEnerg}).

Now we explain the main ideas of this paper.

In \cite{tao}, Tao was able to bound on arbitrary long time intervals the $L_{t}^{4} L_{x}^{12}$ norm of solutions of the energy-critical
equation (\ref{Eqn:QuintDefocWave}) by a quantity that depends exponentially on their energy. This estimate can be viewed as a long time
estimate. Unfortunately we cannot expect to prove a similar result for ( \ref{Eqn:LogQuint} ) since we are not in the energy-critical regime.
However we shall prove the following proposition

\begin{prop}{\textbf{"Long time estimate"}}

Let $J=[t_{1},t_{2}]$ be a time interval. Let $u$ be a \textit{classical solution} of (\ref{Eqn:LogQuint}). Assume that

\begin{equation}
\begin{array}{ll}
\left\| u \right\|_{L_{t}^{\infty} \tilde{H}^{2}( J \times
\mathbb{R}^{3}) } & \leq M
\end{array}
\label{Eqn:ApEstH2norm}
\end{equation}
holds for some $M \geq 0$. Then there exist three constants $C_{L,0}
> 0$, $C_{L,1} > 0$ and $C_{L,2} > 0$ such that

\begin{itemize}

\item if $E << \frac{1}{g^{\frac{1}{2}}(M)}$ (small energy regime) then

\begin{equation}
\begin{array}{ll}
\| u \|^{4}_{L_{t}^{4} L_{x}^{12} (J \times \mathbb{R}^{3})} & \leq
C_{L,0}
\end{array}
\label{Eqn:LonEstSmNrj}
\end{equation}

\item if $E \gtrsim \frac{1}{g^{\frac{1}{2}} (M)}$ (large energy regime) then

\begin{equation}
\begin{array}{ll}
\| u \|^{4}_{L_{t}^{4} L_{x}^{12} (J \times \mathbb{R}^{3})} & \leq
\left( C_{L,1} ( E g(M) ) \right)^{C_{L,2} (E^{\frac{193}{4}+}
g^{\frac{225}{8}+}(M) )}
\end{array}
\label{Eqn:LonEst}
\end{equation}

\end{itemize}
\label{prop:LongEst}
\end{prop}
This proposition shows that we can control the $L_{t}^{4} L_{x}^{12}
(J \times \mathbb{R}^{3})$ norm of  solutions of
(\ref{Eqn:LogQuint}) by their energy and an \textit{a priori} bound
of their $L_{t}^{\infty} \tilde{H}^{2} ( J \times \mathbb{R}^{3})$
norm. We would like to control the pointwise-in-time $\tilde{H}^{2}
(\mathbb{R}^{3}) \times H^{1}(\mathbb{R}^{3})$ norm of $u$ on an
interval $[0,T]$, with $T$ arbitrary large. This is done by an
induction on time. We assume that this norm is controlled on $[0,T]$
by a number $M_{0}$. Then by continuity we can find a slightly
larger interval $[0,T^{'}]$ such that this norm is bounded by (say)
$2M_{0}$ on $[0,T^{'}]$. This is our \textit{a priori} bound. We
subdivide $[0,T^{'}]$ into subintervals where the $L_{t}^{4}
L_{x}^{12}$ norm of $u$ is small and we control the
pointwise-in-time $\tilde{H}^{2} (\mathbb{R}^{3}) \times
H^{1}(\mathbb{R}^{3})$ norm of $u$ on each of these subintervals
(see Lemma \ref{lem:LocalBound}). Since $g$ varies slowly we can
estimate the number of intervals of this partition by using
Proposition \ref{prop:LongEst} and we can prove \textit{a
posteriori} that $ \| u(t) \|_{\tilde{H}^{2}(\mathbb{R}^{3})} + \|
\partial_{t} u(t) \|_{\tilde{H}^{1}(\mathbb{R}^{3})} $ is bounded on $[0,T^{'}]$ by $M_{0}$,
provided that $M_{0}$ is large enough: see
Section \ref{sec:ProofThm}.\\
The proof of Proposition \ref{prop:LongEst} is a modification of the
argument used in \cite{tao} to establish a tower-exponential bound
of the $L_{t}^{4} L_{x}^{12}(J \times \mathbb{R}^{3})$ norm of $v$,
solution of (\ref{Eqn:QuintDefocWave}). We divide $J$ into
subintervals $J_{i}$ where the $L_{t}^{4} L_{x}^{12}$ norm of $u$,
solution of (\ref{Eqn:LogQuint}), is ``substantial''. Then by using
the Strichartz estimates and the Sobolev embedding
(\ref{Eqn:SobEmbed}) we notice that the  $L_{t}^{\infty}
L_{x}^{6}(J_{i} \times \mathbb{R}^{3})$ norm of $u$ is also
``substantial'', more precisely we find a lower bound that depends
on the energy $E$ and $g(M)$. Then by Proposition \ref{prop:InvSob}
we can localize a bubble where the mass concentrates and we prove
that the size of these subintervals is also ``substantially'' large.
Tao \cite{tao} used the mass concentration to construct a solution
$\tilde{v}$ of (\ref{Eqn:QuintDefocWave}) that has a smaller energy
than $v$ and that coincides with $v$ outside a cone. The idea behind
that is to use an induction on the levels of energy, due to Bourgain
\cite{bour}, and the small energy theory following from the
Strichartz estimates in order to control the $L_{t}^{4} L_{x}^{12}$
norm of $v$ outside a cone. Unfortunately it seems almost impossible
to apply this procedure to our problem. Indeed the energy of the
constructed solution $\tilde{u}$ is smaller than the energy $E$ of
$u$ by an amount that depends on $E$ but also on $g(M)$ and
therefore an induction on the levels of the energy is possible if
the $L_{t}^{\infty} \tilde{H}^{2}(J \times \mathbb{R}^{3})$ norm of
$\tilde{u}$ can be controlled by $M$, which is far from being
trivial. It turns out that we do not need to use the Bourgain
induction method. Indeed since we know that the size of the
subintervals $J_{i}$ s is substantially large and since we have a
good control of the $L_{t}^{4} L_{x}^{12}$ norm on these
subintervals it suffices to find an upper bound of the size of their
union in order to conclude. To this end we divide a cone containing
the ball where the mass concentrates and the $J_{i}$ s into
truncated-in-time cones where the $L_{t}^{4} L_{x}^{12}$ norm of $u$
is ``substantial''. Let $\tilde{J}_{1}$,$\tilde{J}_{2}$,... be the
sequence of time intervals resulting from this partition. The mass
concentration helps us to control the size of the first time
interval $\tilde{J}_{1}$. By using an asymptotic stability result we
can prove, roughly speaking, that if we consider two successive
subintervals $\tilde{J}_{j}$, $\tilde{J}_{j+1}$ resulting from this
partition of the cone then the size of $\tilde{J}_{j+1}$ can be
controlled by the size of $\tilde{J}_{j}$: see
(\ref{Eqn:FractTildeJi}). But a potential energy decay estimate
shows that if the size of the union of the $J_{i}$ s is too large
then we can find a large subinterval $[t^{'}_{1}, t^{'}_{2}]$ such
that the $L_{t}^{4} L_{x}^{12}$ norm of $u$ on the cone truncated to
$[t_{1}^{'}, t_{2}^{'}]$ is small. Therefore $[t_{1}^{'}.
t_{2}^{'}]$ cannot be covered by many $\tilde{J}_{j}$ s and one of
them is very large in comparison with its predecessor, which
contradicts (\ref{Eqn:FractTildeJi}). At the end of the process we
can find an upper bound of the size of the union of the subintervals
$J_{i}$ s and consequently we can control the $L_{t}^{4} L_{x}^{12}$
norm of $u$ on the interval $J$.

\begin{rem}
Throughout the paper we frequently use the $x+$ and $x-$ notations.
Indeed the point $(2, \infty)$ is not wave admissible. Therefore we
will work with the point $(2+, \infty -)$: see (\ref{Eqn:Useeps1})
and (\ref{Eqn:Useeps2}). This generates slight variations of many
quantities throughout this paper. Sometimes we might deal with
quantities like $z:=\frac{x+}{y-}$. We cannot conclude directly that
$z=\frac{x}{y}+$. In this case we create a variation of $y$ so small
(comparing to that of $x$) that we have $z=\frac{x}{y}+$. These
details have been omitted for the sake of readability. We strongly
recommend that the reader ignores these slight variations at the
first reading.
\end{rem}

$\textbf{Acknowledgements}:$ The author would like to thank Terence
Tao for suggesting him  this problem.

\section{Proof of Theorem \ref{Thm:GlobLogQuint} }
\label{sec:ProofThm}

The proof relies upon Proposition \ref{prop:LongEst} and the following lemma that we prove in the next subsection.

\begin{lem}{\textbf{"Local boundedness"}}
Let $J=[t_{1},t_{2}]$ be an interval. Assume that $u$ is a
\textit{classical solution} of (\ref{Eqn:LogQuint}). Let $Z(t):= \|
(u(t) , \partial_{t} u(t)) \|_{ \tilde{H}^{2} (\mathbb{R}^{3})
\times H^{1}( \mathbb{R}^{3})} $. There exists $0 < \epsilon << 1$
constant such that if

\begin{equation}
\begin{array}{ll}
\| u \|_{L_{t}^{4} L_{x}^{12} (J \times \mathbb{R}^{3})} & \leq \frac{\epsilon}{ g^{\frac{1}{4}} \left( Z(t_{1})  \right)  }
\end{array}
\label{Eqn:CdtionSmallLt4Lx12}
\end{equation}
then there exists $C_{l} > 0$ such that

\begin{equation}
\begin{array}{ll}
Z(t) & \leq 2 C_{l} Z(t_{1})
\end{array}
\label{Eqn:LocalBound}
\end{equation}
for $t \in J$.
\label{lem:LocalBound}
\end{lem}

We claim that the following set

\begin{equation}
\begin{array}{ll}
\mathcal{F} & : = \left\{ T \in [0, \, \infty ): \,  \sup_{t \in
[0,T]} \| \left( u(t), \partial_{t} u(t)  \right)
\|_{\tilde{H}^{2}(\mathbb{R}^{3}) \times H^{1} (\mathbb{R}^{3})}
\leq M_{0} \right\}
\end{array}
\end{equation}
is equal  to $[0, \, \infty)$ for some constant $M_{0}:=M_{0}(\|
u_{0} \|_{\tilde{H}^{2}(\mathbb{R}^{3})}, \| u_{1}
\|_{H^{1}(\mathbb{R}^{3})})$ large enough. Indeed

\begin{itemize}

\item $ 0 \in \mathcal{F} $: clear

\item $\mathcal{F}$ is closed by continuity

\item $\mathcal{F}$ is open. Indeed let $T \in \mathcal{F}$. Then by continuity there exists $\delta>0$ such that

\begin{equation}
\begin{array}{ll}
\sup_{t \in [0, T^{'}]}  \| \left( u(t), \partial_{t} u(t)  \right) \|_{H^{2}(\mathbb{R}^{3}) \times H^{1} (\mathbb{R}^{3})} & \leq 2 M_{0}
\end{array}
\label{Eqn:InducBdH2}
\end{equation}
for every $T^{'} \in [0,T + \delta)$. By (\ref{Eqn:LonEstSmNrj}) and
(\ref{Eqn:LonEst}) we have

\begin{equation}
\begin{array}{ll}
\| u \|^{4}_{L_{t}^{4} L_{x}^{12} ([0, T^{'}] \times
\mathbb{R}^{3})} & \leq \max{ \left( C_{L,0},  ( C_{L,1} E \,
g(2M_{0}))^{C_{L,2} (E^{\frac{193}{4}+} g^{\frac{225}{8}+}(2 M_{0})
)} \right) }
\end{array}
\label{Eqn:AprLongEst}
\end{equation}
Let $N \geq 1$ and let $\underline{Z}(0):= \max{ \left( Z(0), \, 1
\right)}$. Without loss of generality we can assume that $C_{l}
>> 1$ so that $2 C_{l} \underline{Z}(0) >> 1$ and $\log^{c} \left( 2 C_{l} \underline{Z}(0) \right) >>
1$. We have, by the elementary rules of the logarithm \footnote{
such as the product rule $\log{(a_{1} a_{2}) = \log{(a_{1})}} +
\log{(a_{2})}$} and the inequality  $\log^{c}(2nx) \leq \log^{c}
\left( (2n)^{x}   \right)$ for $n \geq 1$ and $x >> 1$

\begin{equation}
\begin{array}{ll}
\sum \limits_{n=1}^{N} \frac{\epsilon^{4}}{g \left(
(2C_{l})^{n}Z_{0} \right)} & \geq \sum \limits_{n=1}^{N}
\frac{\epsilon^{4}}
{\log^{c} \left( \log ( (2C_{l})^{2n} \underline{Z}^{2n}(0) +10 ) \right) } \\
& \gtrsim \sum \limits_{n=1}^{N}  \frac{1}{\log^{c} \left( 2n \log{( 2 C_{l} \underline{Z}(0))}   \right) } \\
& \gtrsim \frac{1}{\log^{c} \left( 2 C_{l} \underline{Z}(0) \right)} \sum \limits_{n=1}^{N} \frac{1}{\log^{c}(2n)} \\
& \gtrsim \frac{1}{\log^{c} \left( 2 C_{l} \underline{Z}(0) \right)} \int_{1}^{N+1} \frac{1}{\log^{c}(2t)} \, dt \\
& \gtrsim \frac{1}{\log^{c} \left( 2 C_{l} \underline{Z}(0) \right)}
\int_{1}^{N+1} \frac{1}{t^{\frac{1}{2}}} \, dt \\
& \gtrsim \frac{N^{\frac{1}{2}}}{\log^{c} \left( 2 C_{l}
\underline{Z}(0) \right)}
\end{array}
\label{Eqn:LowerBoundSum}
\end{equation}
By Lemma \ref{lem:LocalBound}, (\ref{Eqn:AprLongEst}) and
(\ref{Eqn:LowerBoundSum}) we can construct a partition $(J_{n})_{1
\leq n \leq N}$ of $[0,T^{'}]$ such that $\| u \|_{L_{t}^{4}
L_{x}^{12} (J_{n} \times \mathbb{R}^{3})} =
\frac{\epsilon}{g^{\frac{1}{4}}\left( (2C_{l})^{n} Z_{0} \right)}$,
$1 \leq n < N $, $\| u \|_{L_{t}^{4} L_{x}^{12} (J_{N} \times
\mathbb{R}^{3})} \leq \frac{\epsilon}{g^{\frac{1}{4}}\left(
(2C_{l})^{N} Z_{0} \right)}$, $Z(t) \leq (2C_{l})^{n} Z(0)$ for $t
\in J_{1} \cup... \cup J_{n}$ and

\begin{equation}
\begin{array}{ll}
\frac{N^{\frac{1}{2}}}{\log^{c} \left( 2 C_{l} \underline{Z}(0)
\right)} & \leq \max{ \left( C_{L,0}, ( C_{L,1} E \,
g(2M_{0}))^{C_{L,2} (E^{\frac{193}{4}+} g^{\frac{225}{8}+}(2 M_{0})
)} \right) }
\end{array}
\label{Eqn:ConstrN}
\end{equation}
Since $c < \frac{8}{225}$ we have by (\ref{Eqn:BoundE})

\begin{equation}
\begin{array}{ll}
\log{N} & \lesssim  \log^{c} \left( 2 C_{l} \underline{Z}(0) \right)
+  \log{(C_{L,0})} +  C_{L,2} E^{\frac{193}{4}+}
\log^{\frac{225c}{8}+} \log{(10 + 4 M^{2}_{0})}
\log{ \left( C_{L,1}  E  \log^{c} \log  (10 + 4 M_{0}^{2}) \right) } \\
& \leq  \log \left( \frac{ \log{ \left( \frac{M_{0}}{Z(0)}
\right)}}{\log{(2 C_{l})}} \right)
\end{array}
\label{Eqn:ControlNorm}
\end{equation}
if $M_{0}=M_{0} \left( \| u_{0} \|_{\tilde{H}^{2}(\mathbb{R}^{3})},
\| u_{1} \|_{H^{1}(\mathbb{R})} \right)$ is large enough. To prove
the last inequality in (\ref{Eqn:ControlNorm}) it is enough, by
using (\ref{Eqn:BoundE}), to notice that $ \lim_{M_{0} \rightarrow
\infty} f(M_{0}) =0$ with

\begin{equation}
\begin{array}{ll}
f(M_{0}) & : = \frac {\log^{c} \left( 2 C_{l} \underline{Z}(0)
\right) +  \log{(C_{L,0})} +  C_{L,2} E^{\frac{193}{4}+}
\log^{\frac{225c}{8}+} \log{(10 + 4 M^{2}_{0})} \log{ \left( C_{L,1}
E  \log^{c} \log  (10 + 4 M_{0}^{2}) \right) }} { \log \left( \frac{
\log{ \left( \frac{M_{0}}{Z(0)} \right)}}{\log{(2 C_{l})}} \right)}
\end{array}
\end{equation}

Therefore we conclude that

\begin{equation}
\begin{array}{ll}
\sup_{t \in [0,T^{'}]} \left\| (u(t), \partial_{t}u(t))  \right\|_{H^{2}(\mathbb{R}^{3}) \times H^{1}(\mathbb{R}^{3})} & \leq (2 C_{l})^{N} Z(0) \\
& \leq M_{0}
\end{array}
\end{equation}

\end{itemize}

\subsection{Proof of Lemma \ref{lem:LocalBound}}

By the Strichartz estimates (\ref{Eqn:Str}), the Sobolev embeddings
(\ref{Eqn:SobEmbed}) and (\ref{Eqn:SobEmbed2})  and the elementary
estimate $|u^{5} \nabla \left( g(u) \right)| \lesssim |u^{4} \nabla
u g(u)|$

\begin{equation}
\begin{array}{ll}
Z(t) & \lesssim Z(t_{1}) + \| u^{5} g(u) \|_{L_{t}^{1}L_{x}^{2} ( [t_{1},t] \times \mathbb{R}^{3} )} + \| u^{4}  \nabla u g(u) \|_{L_{t}^{1}
L_{x}^{2} ([t_{1},t] \times \mathbb{R}^{3})} + \| u^{5} \nabla ( g(u)) \|_{L_{t}^{1} L_{x}^{2} ( [t_{1},t] \times \mathbb{R}^{3})} \\ \\
&
\lesssim Z(t_{1}) + \| u^{5} g(u) \|_{L_{t}^{1}L_{x}^{2} ( [t_{1},t] \times \mathbb{R}^{3} )} + \| u^{4}  \nabla u g(u) \|_{L_{t}^{1}
L_{x}^{2} ( [t_{1},t] \times \mathbb{R}^{3})} \\ \\
& \lesssim Z(t_{1}) + \| u \|^{4}_{L_{t}^{4} L_{x}^{12} ( [t_{1},t] \times \mathbb{R}^{3})} \| u \|_{L_{t}^{\infty} L_{x}^{6} ( [t_{1},t] \times
\mathbb{R}^{3})} g \left( \| u \|_{L_{t}^{\infty} L_{x}^{\infty} ( [t_{1},t] \times \mathbb{R}^{3})} \right) \\  & + \| u \|^{4}_{L_{t}^{4}
L_{x}^{12} ( [t_{1},t] \times \mathbb{R}^{3})} \| \nabla u \|_{L_{t}^{\infty} L_{x}^{6} ( [t_{1},t] \times
\mathbb{R}^{3})} g \left( \| u \|_{L_{t}^{\infty} L_{x}^{\infty} ( [t_{1},t] \times \mathbb{R}^{3})} \right) \\ \\
& \lesssim Z(t_{1}) + \| u \|^{4}_{L_{t}^{4} L_{x}^{12} ( [t_{1},t]
\times \mathbb{R}^{3})} Z(t) g ( Z(t)) \\
\end{array}
\label{Eqn:IneqZtImpl}
\end{equation}
Let $C_{l}$ be the constant determined by the last inequality in
(\ref{Eqn:IneqZtImpl}). From (\ref{Eqn:CdtionSmallLt4Lx12}),
(\ref{Eqn:IneqZtImpl}) and a continuity argument we have
(\ref{Eqn:LocalBound}).

\section{Proof of Proposition \ref{prop:LongEst}}

The proof relies upon five lemmas that are proved in the next
sections.

\begin{lem}{\textbf{"Long time estimate if energy small"}}
Let $J=[t_{1},t_{2}]$ be a time interval. Let $u$ be a \textit{classical solution} of (\ref{Eqn:LogQuint}). Assume that (\ref{Eqn:ApEstH2norm})
holds. If

\begin{equation}
\begin{array}{ll}
E & << \frac{1}{g^{\frac{1}{2}}(M)}
\end{array}
\label{Eqn:BoundNrjM}
\end{equation}
then

\begin{equation}
\begin{array}{ll}
\| u \|_{L_{t}^{4} L_{x}^{12} (J \times \mathbb{R}^{3})} & \lesssim 1
\end{array}
\label{Eqn:Lt4Lx12BdNrjSmall}
\end{equation}
\label{lem:Lt4Lx12Nrjsmall}
\end{lem}

\begin{lem}{\textbf{" If $ \| u \|_{L_{t}^{4} L_{x}^{12} (J \times \mathbb{R}^{3})}$ non negligeable then  existence of a mass concentration bubble
and size of $J$ bounded from below"}} Let $u$ be a \textit{classical solution} of (\ref{Eqn:LogQuint}). Let $J$ be a time interval. Assume that
(\ref{Eqn:ApEstH2norm}) holds. Let $\eta$ be a positive number such that

\begin{equation}
\begin{array}{ll}
\eta & \leq \frac{E^{\frac{1}{12}}}{\gfivetwfourM}
\end{array}
\label{Eqn:UpperBdeta}
\end{equation}
If $ \| u \|_{L_{t}^{4} L_{x}^{12} (J \times \mathbb{R}^{3})} \geq \eta$ then

\begin{equation}
\begin{array}{ll}
\| u \|_{L_{t}^{\infty} L_{x}^{6}(J \times \mathbb{R}^{3})} &
\gtrsim \eta^{2+} E^{- \left( \frac{1}{2} + \right)}
\end{array}
\label{Eqn:logthfivetwfourM}
\end{equation}
Moreover there exist a point $x_{0} \in \mathbb{R}^{3}$, a time $t_{0} \in J$ and a positive number $ r $ such that we have the mass
concentration estimate in the ball $B(x_{0},r)$

\begin{equation}
\begin{array}{ll}
\int_{B(x_{0},r)} \left| u(t_{0},y) \right|^{2} \, dy & \gtrsim
\eta^{6+} E^{-(2+)} r^{2}
\end{array}
\label{Eqn:Massconc}
\end{equation}
and the following lower bound on the size of $J$

\begin{equation}
\begin{array}{ll}
|J| & \gtrsim \eta^{4}  E^{-\frac{2}{3}} r
\end{array}
\label{Eqn:LowerBoundSizeJ}
\end{equation}
\label{lem:Lt4Lx12Mass}
\end{lem}

\begin{lem}{\textbf{"Potential energy decay estimate"}}
Let $u$ be a \textit{classical solution} of (\ref{Eqn:LogQuint}).
Let $[a,b]$ be an interval. Then we have the potential energy decay
estimate
\begin{equation}
\begin{array}{ll}
\int_{|x| \leq b} F(u(b,x)) \,dx & \lesssim \frac{a}{b} \left( e(a) + e^{\frac{1}{3}}(a) \right) + e(b)-e(a) + \left( e(b) -e(a)
\right)^{\frac{1}{3}}
\end{array}
\label{Eqn:DecayPot}
\end{equation}
\label{lem:DecayPotNrjCone}
\end{lem}

\begin{lem}{\textbf{" $ L_{t}^{4} L_{x}^{12} $ norm of $u$ is small on a large truncation of the forward light cone"}}
Let $J=[t_{1}, \, t_{2}]$ be an interval. Let $u$ be a \textit{classical solution} of (\ref{Eqn:LogQuint}). Assume that (\ref{Eqn:ApEstH2norm})
holds. Let $\eta$ be a positive number such that

\begin{equation}
\begin{array}{ll}
\eta & << \min \left( E^{\frac{1}{4}},E^{\frac{5}{18}}, \frac{E^{\frac{1}{12}}}{g^{\frac{5}{24}}(M)} \right)
\end{array}
\label{Eqn:EtaHyp}
\end{equation}
Assume also that there exists $ C_{2} >> 1$ such that

\begin{equation}
\begin{array}{ll}
\left[ t_{1}, \, ( C_{2} E^{10+} \eta^{-(36+)}  )^{4 C_{2} E^{10 +}
\eta^{-(36+)}} t_{1} \right] & \subset J
\end{array}
\label{Eqn:ConstSizeSubInt}
\end{equation}
Then there exists a subinterval $J^{'}=[t_{1}^{'}, t_{2}^{'}]$ such
that $ \left| \frac{t_{2}^{'}}{t_{1}^{'}} \right| \sim E^{10+}
\eta^{-(36+)}$ and

\begin{equation}
\begin{array}{ll}
\| u  \|_{L_{t}^{4} L_{x}^{12} (\Gamma_{+}(J^{'}))} & \leq \eta
\end{array}
\end{equation}
\label{lem:Lt4Lx12}
\end{lem}

\begin{lem}{\textbf{"Asymptotic stability"}}
Let $J=[t_{1},t_{2}]$ be a time interval. Let $J^{'}=[t^{'}_{1}, \, t^{'}_{2}] \subset J$ and let $t \in J/J^{'}$. Let $u$ be a \textit{
classical solution} of (\ref{Eqn:LogQuint}). Assume that (\ref{Eqn:ApEstH2norm}) holds. Then

\begin{equation}
\begin{array}{ll}
\| u_{l,t^{'}_{2}}(t) - u_{l,t^{'}_{1}}(t) \|_{L_{x}^{\infty} (\mathbb{R}^{3})} & \lesssim \frac{E^{\frac{5}{6}}
\gonesixM}{dist^{\frac{1}{2}}(t,J^{'})}
\end{array}
\label{Eqn:AsympStab}
\end{equation}
\label{lem:AsympStab}
\end{lem}

We are ready to prove Proposition \ref{prop:LongEst}. We assume that
we have an \textit{a priori} bound $M$ of the $L_{t}^{\infty}
\tilde{H}^{2} (J \times \mathbb{R}^{3})$ norm  of the solution $u$.
There are two steps

\begin{itemize}

\item If $E << \frac{1}{g^{\frac{1}{2}}(M)}$ then we know from Lemma \ref{lem:Lt4Lx12Nrjsmall} that
(\ref{Eqn:LonEstSmNrj}) holds.

\item Therefore we assume that the energy is large, i.e

\begin{equation}
\begin{array}{ll}
E & \gtrsim \frac{1}{g^{\frac{1}{2}}(M)}
\end{array}
\label{Eqn:LowerBoundE}
\end{equation}
We can assume without loss of generality that

\begin{equation}
\begin{array}{ll}
\| u \|_{L_{t}^{4}  L_{x}^{12} (J \times \mathbb{R}^{3})} & \geq \frac{E^{\frac{1}{12}}}{g^{\frac{5}{24}}(M)}
\end{array}
\label{Eqn:BdAssLt4Lx12}
\end{equation}
From (\ref{Eqn:BdAssLt4Lx12}) we can partition $J$ into subintervals $J_{1}$, ..., $J_{l}$ such that for $i=1, ..., l-1$

\begin{equation}
\begin{array}{ll}
\| u \|_{L_{t}^{4} L_{x}^{12} (J_{i} \times \mathbb{R}^{3})} & = \frac{E^{\frac{1}{12}}}{g^{\frac{5}{24}}(M)}
\end{array}
\label{Eqn:AssJieq1}
\end{equation}
and

\begin{equation}
\begin{array}{ll}
\| u \|_{L_{t}^{4}  L_{x}^{12} (J_{l} \times \mathbb{R}^{3})} & \leq \frac{E^{\frac{1}{12}}}{g^{\frac{5}{24}}(M)}
\end{array}
\label{Eqn:AssJieq2}
\end{equation}
Before moving forward we say that an interval $J_{i}$ is
\textit{exceptional} if \footnote{The numbers $\frac{193}{4}$ and
$\frac{225}{8}$ will play an important role in
(\ref{Eqn:ConditionT2T1}) }

\begin{equation}
\begin{array}{ll}
\| u_{l,t_{1}} \|_{L_{t}^{4} L_{x}^{12} (J_{i} \times
\mathbb{R}^{3})} + \| u_{l,t_{2}} \|_{L_{t}^{4} L_{x}^{12} (J_{i}
\times \mathbb{R}^{3})} & \geq \frac{1}{ \left( C_{3} E g(M)
\right)^{C_{4} \left( E^{\frac{193}{4} +} g^{\frac{225}{8} +}(M)
\right)}}
\end{array}
\label{Eqn:DefExcep}
\end{equation}
for some $C_{3} >> 1$, $C_{4} >> 1$ to be chosen later. Otherwise $J_{i}$ is  \textit{unexceptional}. Let $\mathcal{E}$ denote the set of
$J_{i}^{'}$ s that are exceptional and let $\overline{ \mathcal{E}^{c}}$ denote the set of nonempty sequences of consecutive unexceptional
intervals $J_{i}$. By (\ref{Eqn:Str}), (\ref{Eqn:LowerBoundE}) and (\ref{Eqn:DefExcep})

\begin{equation}
\begin{array}{ll}
\card{(\mathcal{E})}  & \lesssim E^{2}  \left[ O (  E g(M) ) \right]^{O \left( E^{\frac{193}{4} +} g^{\frac{225}{8} +}(M)  \right)} \\
& \lesssim \left[ O (  E g(M) ) \right]^{O \left( E^{\frac{193}{4}
+} g^{\frac{225}{8} +}(M)  \right)}
\end{array}
\end{equation}
Since $\card{(\overline{\mathcal{E}^{c}})} \lesssim \card{ \left( \mathcal{E} \right) }$ we have

\begin{equation}
\begin{array}{ll}
\| u \|^{4}_{L_{t}^{4} L_{x}^{12}(J \times \mathbb{R}^{3})} &
\lesssim \left[ O (E g(M)) \right]^{O \left( E^{\frac{193}{4}+}
g^{\frac{225}{8} +}(M) \right)}  \left(
\frac{E^{\frac{1}{3}}}{g^{\frac{5}{6}}(M)} + \sup_{K \in \overline
{\mathcal{E}^{c}}} \| u \|^{4}_{L_{t}^{4} L_{x}^{12} (K \times
\mathbb{R}^{3})} \right)
\end{array}
\label{Eqn:ExUnexuLt4Lx12}
\end{equation}
Let $ K=J_{i_{0}} \cup ...\cup J_{i_{1}}$ is a sequence of
consecutive unexceptional intervals. If $N(K)$ is the number of
$J_{i}$ s making $K$ then by (\ref{Eqn:LowerBoundE}),
(\ref{Eqn:AssJieq1}), (\ref{Eqn:AssJieq2}) and
(\ref{Eqn:ExUnexuLt4Lx12}) we have

\begin{equation}
\begin{array}{ll}
\| u \|_{L_{t}^{4} L_{x}^{12}(J \times \mathbb{R}^{3})} & \lesssim
\left( \sup_{K \in \overline {\mathcal{E}^{c}} } N(K) \right) \left[
O (E g(M)) \right]^{O \left( E^{\frac{193}{4} +} g^{\frac{225}{8}
+}(M) \right)}
\end{array}
\label{Eqn:CrudeEstLt4Lx12}
\end{equation}
Therefore it suffices to estimate $N(K)$ for every $K=J_{i_{0}} \cup
.... \cup J_{i_{1}}$. We will do that by first determining a lower
bound for the size of the elements $J_{i}$' s and then by estimating
the size of $K$.  By (\ref{Eqn:LowerBoundE}) ,(\ref{Eqn:AssJieq1}),
(\ref{Eqn:AssJieq2}) and Lemma \ref{lem:Lt4Lx12Mass} there exists
for $i \in [i_{0}, \, ... i_{1}]$ a $(t_{i}, r_{i}, x_{i}) \in
\left( J_{i} \times ( \, 0, \infty) \times \mathbb{R}^{3} \right)$
such that

\begin{equation}
\begin{array}{ll}
\frac{1}{r^{2}_{i}} \int_{B(x_{i},r_{i})} |u(t_{i},y)|^{2} \, dy &
\gtrsim \frac{E^{-\left( \frac{3}{2} + \right)}}{g^{ \frac{5}{4}+}
(M)}
\end{array}
\label{Eqn:massinfIi}
\end{equation}
and

\begin{equation}
\begin{array}{ll}
|J_{i}| & \gtrsim  \frac{E^{-\frac{1}{3}} r_{i}} {g^{\frac{5}{6}}(M)}
\end{array}
\end{equation}

Let $k \in [i_{0}, \ .. i_{1}]$ such that $r_{k} = \min_{_{i \in [i_{0}, \, i_{1}]}} r_{i}$, let $f(t,r,x):=\frac{1}{r^{2}} \int_{B(x,r)}
|u(t,y)|^{2} \, dy$, let $C_{5}$ be the constant determined by (\ref{Eqn:massinfIi}). Let $r_{0}=r_{0}(M)$ such that \\
$r_{0} M^{2} = \frac{C_{5} E^{- \left( \frac{3}{2} + \right)}}{4 g^{
\frac{5}{4} + }(M)} $. Since $f(t,r,x) \leq r M^{2}$ we have $
f(t,r_{0},x) \leq \frac{C_{5} E^{- \left( \frac{3}{2} + \right)}}{4
g^{\frac{5}{4}+}(M)} $. The set $A:= \left\{ (t,r,x): \, t \in K, \,
r_{0} \leq r \leq r_{k}, \, x \in \mathbb{R}^{3} \right\}$ is
connected. Therefore its image is connected by $f$ and there exists
$\left( \tilde{t}, \, \tilde{r}, \, \tilde{x} \right) \in K \times
[r_{0}, \, r_{k}] \times \mathbb{R}^{3}$ such that $ f(\tilde{t},
\tilde{r}, \tilde{x}) = \frac{C_{5} E^{- \left( \frac{3}{2} +
\right) } }{2 g^{\frac{5}{4} +}(M)}$ . In other words we have the
following mass concentration

\begin{equation}
\begin{array}{ll}
\frac{1}{\tilde{r}^{2}} \int_{B(\tilde{x},\tilde{r})}
u^{2}(\tilde{t}, x) \, dx & =  \frac{C_{5} E^{- \left( \frac{3}{2} +
\right) }}{2 g^{\frac{5}{4} +}(M)}
\end{array}
\label{Eqn:MassconcPr}
\end{equation}
Moreover we have the useful lower bound for the size of $J_{i}$, $
i_{0} \leq i \leq i_{1}$  \footnote{ Notice that we have the  lower
bound $\tilde{r} \geq \frac{C_{5} E^{-\left( \frac{3}{2}+
\right)}}{4 M^{2} g^{\frac{5}{4}+}(M)}$. One might think that the
presence of $\tilde{r}$ in (\ref{Eqn:EstSizeJj}) is annoying since
this lower bound is crude. However we will see that $\tilde{r}$
disappears at the end of the process: see (\ref{Eqn:EstNK}).
Therefore a sharp lower bound is not required.}

\begin{equation}
\begin{array}{ll}
|J_{i}| & \gtrsim \tilde{r}
\frac{E^{-\frac{1}{3}}}{g^{\frac{5}{6}}(M)}
\end{array}
\label{Eqn:EstSizeJj}
\end{equation}
At this point we need to use the following lemma that we will prove
in the next subsection

\begin{lem}
Let K be a sequence of unexceptional intervals. Assume that there exists $(\bar{t}, \, \bar{x}, \bar{r}) \in K \times \mathbb{R}^{3} \times (0,
\infty)$ such that

\begin{equation}
\begin{array}{ll}
\frac{1}{\bar{r}^{2}} \int_{B(\bar{x},\bar{r})} u^{2}(\bar{t},y) \,
dy & \gtrsim E^{- \left( \frac{3}{2} + \right)}{g^{\frac{5}{4} +}
(M)}
\end{array}
\end{equation}
Then there exist two constants $C_{6} >> 1$, $C_{7}>>1$ such that

\begin{equation}
\begin{array}{ll}
|K| & \leq  \left( C_{6} E g(M) \right)^{C_{7} E^{ \frac{193}{4} +}
g^{\frac{225}{8} +}(M) } \bar{r}
\end{array}
\label{Eqn:EstSizeK}
\end{equation}
\label{Lem:EstDiamK}
\end{lem}
Lemma \ref{Lem:EstDiamK} gives information about the size of $K$. We
postpone the proof of this lemma to the next subsection. If we
combine it to (\ref{Eqn:EstSizeJj}) we can estimate $N(K)$. More
precisely by Lemma \ref{Lem:EstDiamK}, (\ref{Eqn:EstSizeJj}) and
(\ref{Eqn:LowerBoundE}) we have

\begin{equation}
\begin{array}{ll}
N(K) & \lesssim \frac{ \left( C_{6} E g(M) \right)^{C_{7}
E^{\frac{193}{4}+} g^{\frac{225}{8}+}(M)} \tilde{r}} {\tilde{r}
\frac{E^{-\frac{1}{3}}}{g^{\frac{5}{6}}(M)}} \\
& \lesssim \left( O( E g(M)) \right)^{O( E^{\frac{193}{4}+}
g^{\frac{225}{8}+}(M))}
\end{array}
\label{Eqn:EstNK}
\end{equation}
Plugging this upper bound of $N(K)$ into (\ref{Eqn:CrudeEstLt4Lx12}) we get (\ref{Eqn:LonEst}).

\subsection{Proof of Lemma \ref{Lem:EstDiamK}}

By using the space translation invariance of (\ref{Eqn:LogQuint}) we
can reduce to the case $\bar{x}=0$ \footnote{we consider the
function $u_{1}(t,x)=u(t,x-\bar{x})$ and we abuse notation in the
sequel by writing $u_{1}$ for $u$ }. By using the time reversal
invariance and the time translation invariance \footnote{we consider
the function $u_{2}(t,x):=u(2 \bar{t}-t,x)$ and we abuse notation in
the sequel by writing $u_{2}$ for $u$ } it suffices to estimate $| K
\cap  [ \bar{t}, \, \infty) |$. By using the time translation
invariance again \footnote{we consider the function
$u_{3}(t,x):=u(t+ (\bar{t}-\bar{r}),x)$ and we abuse notation in the
sequel by writing $u_{3}$ for $u$ } we can assume that
$\bar{t}=\bar{r}$ and therefore $\bar{r} \in K$. Let $ K_{+} :=  K
\cap [\bar{r}, \infty) $. We are interested in estimating $|K_{+}|$.
We would like to use Lemma \ref{lem:Lt4Lx12}. Therefore, we consider
the set $\Gamma_{+}(K_{+})$. We have

\begin{equation}
\begin{array}{ll}
\frac{1}{\bar{r}^{2}} \int_{B(0,\bar{r})} |u(\bar{r},y)|^{2} \, dy &
\gtrsim \frac{E^{- \left(\frac{3}{2} + \right)}
}{g^{\frac{5}{4}+}(M)}
\end{array}
\label{Eqn:ConcMassR}
\end{equation}
Therefore by Proposition \ref{prop:tao} and (\ref{Eqn:ConcMassR}) we have

\begin{equation}
\begin{array}{ll}
\int_{B(0,\bar{r})} |u (t,y)|^{2} \, dy & \gtrsim \frac{E^{-\left(
\frac{3}{2} + \right)} \bar{r}^{2} }{g^{\frac{5}{4} +}(M)}
\end{array}
\end{equation}
if $ (t-\bar{r}) E^{\frac{1}{2}} \leq \frac{c_{0} E^{- \left(
\frac{3}{4} + \right)} \bar{r}}{g^{\frac{5}{8} +}(M)}$ for some
$c_{0} << 1$. Therefore by H\"older there exists $0 < c_{1} << 1$
small enough such that

\begin{equation}
\begin{array}{ll}
\| u \|_{L_{t}^{4} L_{x}^{12} \left( \Gamma_{+} \left(
\left[\bar{r}, \bar{r} + \frac{c_{0} E^{- \left( \frac{5}{4} +
\right)} \bar{r}}{g^{\frac{5}{8} +}(M)} \right] \right) \right) } &
\geq c_{1} \frac{E^{-\frac{17}{16}}}{g^{\frac{25}{32}}(M)}
\end{array}
\label{Eqn:Lt4Lx12RInt}
\end{equation}

\begin{itemize}

\item  $\| u \|_{L_{t}^{4} L_{x}^{12} (\Gamma_{+} (K_{+}))} \leq c_{1} \frac{E^{-\frac{17}{16}}}{g^{\frac{25}{32}}(M)}$. In this case we get
from ( \ref{Eqn:Lt4Lx12RInt} )

\begin{equation}
\begin{array}{ll}
K_{+} \subset \left[ \bar{r}, \, \bar{r} + \frac{c_{0} E^{- \left(
\frac{5}{4} + \right)} \bar{r}}{g^{\frac{5}{8} +}(M)} \right]
\end{array}
\end{equation}
and, using also (\ref{Eqn:LowerBoundE}), we get
(\ref{Eqn:EstSizeK}).

\item $\| u \|_{L_{t}^{4} L_{x}^{12} (\Gamma_{+} (K_{+}))} \geq c_{1} \frac{E^{- \left( \frac{17}{16} + \right)}}{g^{\frac{25}{32} +}(M)}$. In the sequel we denote
by $\tilde{\eta}$ the following number

\begin{equation}
\begin{array}{ll}
\tilde{\eta} := \frac{c_{1}}{4} \frac{E^{- \left( \frac{17}{16}
\right)+}}{g^{\frac{25}{32}+}(M)}
\end{array}
\label{Eqn:ValueTildeEta}
\end{equation}
and we divide $\Gamma_{+}(K_{+})$ into consecutive cone truncations $\Gamma_{+}(\tilde{J}_{1})$,... $\Gamma_{+}(\tilde{J}_{k})$ such that for
$j=1,..,k-1$

\end{itemize}

\begin{equation}
\begin{array}{lll}
\| u \|_{L_{t}^{4} L_{x}^{12} ( \Gamma_{+} ( \tilde{J}_{j}) )} & = \tilde{\eta}
\end{array}
\label{Eqn:CdtildeJi}
\end{equation}
and

\begin{equation}
\begin{array}{ll}
\| u \|_{L_{t}^{4} L_{x}^{12} (\Gamma_{+} ( \tilde {J}_{k}))} & \leq \tilde{\eta}
\end{array}
\label{Eqn:CdtildeJk}
\end{equation}
We get from (\ref{Eqn:Lt4Lx12RInt})

\begin{equation}
\begin{array}{ll}
\tilde{J}_{1} & \subset \left[ \bar{r}, \bar{r} + \frac{c_{0} E^{-
\left( \frac{5}{4} + \right)} \bar{r}}{g^{\frac{5}{8}+}(M)} \right]
\end{array}
\label{Eqn:BdFstInt}
\end{equation}
Now we prove the following result before moving forward

\begin{res}
If $j \in [1,...,k-1]$ we either have

\begin{equation}
\begin{array}{ll}
|\tilde{J}_{j+1}| & \lesssim |\tilde{J}_{j}|  \tilde{\eta}^{-4} E^{\frac{8}{3}} g^{\frac{1}{3}}(M)
\end{array}
\label{Eqn:FractTildeJi}
\end{equation}
or

\begin{equation}
\begin{array}{ll}
|\tilde{J}_{j}| & \geq   \left( C_{6} E g(M) \right)^{C_{7}
E^{\frac{193}{4}+} g^{\frac{225}{8}+}(M)} \bar{r}
\end{array}
\label{Eqn:SizeTildeJi}
\end{equation}
for some  constants $C_{6} >> 1$, $C_{7} >>1$. \label{Res:LacEst}
\end{res}

\begin{proof}
We get from (\ref{Eqn:Strcone2}), (\ref{Eqn:LowerBoundE}) and
(\ref{Eqn:ValueTildeEta})

\begin{equation}
\begin{array}{ll}
\| u - u_{l,t_{j+1}} \|_{L_{t}^{4} L_{x}^{12} (\Gamma_{+} (\tilde{J}_{j}))} & \lesssim \| u^{5} g(u) \|_{L_{t}^{1} L_{x}^{2}
(\Gamma_{+}(\tilde{J_{j}} \cup \tilde{J_{j+1}} ) ) } \\
& \lesssim \| u^{4} \|_{L_{t}^{1} L_{x}^{3} (\Gamma_{+}(\tilde{J_{j}} \cup \tilde{J_{j+1}} ) ) }
 \| u g^{\frac{1}{6}}(u) \|_{L_{t}^{\infty} L_{x}^{6} (\Gamma_{+} (\tilde{J_{j}} \cup \tilde{J_{j+1}} ))}
g^{\frac{5}{6}}(M) \\
& \lesssim \tilde{\eta}^{4} E^{\frac{1}{6}} g^{\frac{5}{6}}(M) \\
& << \tilde{\eta}
\end{array}
\end{equation}
with $J_{j}=[t_{j-1},t_{j}]$. Therefore  by (\ref{Eqn:CdtildeJi}) we
have $ \| u_{l,t_{j+1}} \|_{L_{t}^{4} L_{x}^{12}
(\Gamma_{+}(\tilde{J}_{j}))} \sim \tilde{\eta}$. This implies that

\begin{equation}
\begin{array}{ll}
\| u_{l,t_{j+1}} - u_{l,t_{2}} \|_{L_{t}^{4} L_{x}^{12} (\Gamma_{+}(\tilde{J}_{j}))} & \gtrsim \tilde{\eta}.
\end{array}
\label{Eqn:Diffutit2}
\end{equation}
or

\begin{equation}
\begin{array}{ll}
\| u_{l,t_{2}} \|_{L_{t}^{4} L_{x}^{12} (\Gamma_{+}(\tilde{J}_{j}))} & \gtrsim \tilde{\eta}
\end{array}
\label{Boundut2}
\end{equation}

\begin{itemize}

\item \textit{Case 1}: $\| u_{l,t_{j+1}} - u_{l,t_{2}} \|_{L_{t}^{4} L_{x}^{12} (\Gamma_{+} (\tilde{J}_{i}))} \gtrsim \tilde{\eta} $. By
Lemma \ref{lem:AsympStab} and H\"older we have

\begin{equation}
\begin{array}{ll}
\| u_{l,t_{j+1}} - u_{l,t_{2}} \|_{ L_{t}^{4} L_{x}^{12} (\Gamma_{+} (\tilde{J}_{j}) )} & \lesssim  |\tilde{J}_{j}|^{\frac{1}{4}} \|
u_{l,t_{j+1}}
- u_{l,t_{2}} \|_{L_{t}^{\infty} L_{x}^{12} ( \Gamma_{+} (\tilde{J}_{j}))} \\
& \lesssim | \tilde{J}_{j} |^{\frac{1}{4}} \| u_{l,t_{j+1}} -u_{l,t_{2}} \|^{\frac{1}{2}}_{L_{t}^{\infty} L_{x}^{\infty} (\Gamma_{+}
(\tilde{J}_{j}))} \| u_{l,t_{j+1}} - u_{l,t_{2}} \|^{\frac{1}{2}}_{L_{t}^{\infty} L_{x}^{6} (\Gamma_{+}(\tilde{J}_{j}))} \\
& \lesssim \frac{|\tilde{J}_{j}|^{\frac{1}{4}} E^{\frac{2}{3}}  g^{\frac{1}{12}}(M) } {| \tilde{J}_{j+1} |^{\frac{1}{4}}}
\end{array}
\label{Eqn:Bounddiffut}
\end{equation}
We get (\ref{Eqn:FractTildeJi}) from  (\ref{Eqn:Diffutit2}) and (\ref{Eqn:Bounddiffut}).

\item \textit{Case 2}: $\| u_{l,t_{2}} \|_{L_{t}^{4} L_{x}^{12} (\Gamma_{+} (\tilde{J}_{j})) } \gtrsim \tilde{\eta}$. In this case
$\| u_{l,t_{2}} \|_{L_{t}^{4} L_{x}^{12} (\tilde{J}_{j})} \gtrsim
\tilde{\eta}$. Recall that $K_{+}$ is a subinterval of $K =J_{i_{0}}
\cup.... \cup J_{i_{1}}$, sequence of unexceptional intervals
$J_{i}$, $ i_{0} \leq i \leq i_{1}$. Consequently there are at least
$ \sim \tilde{\eta} \left( C_{3} E g(M) \right)^{ C_{4}
E^{\frac{193}{4} +} g^{\frac{225}{8} +}(M) } $ intervals $J_{j}$
that cover $\tilde{J}_{i}$. Therefore we get (\ref{Eqn:SizeTildeJi})
from (\ref{Eqn:EstSizeJj}) and (\ref{Eqn:LowerBoundE}).

\end{itemize}

\end{proof}

Using Result \ref{Res:LacEst} and Lemma \ref{lem:Lt4Lx12} we can get an upper bound on the size $|K_{+}|$

\begin{res}
We have
\begin{equation}
\begin{array}{ll}
|K_{+}| & \leq  \left( C_{6} E  g(M)  \right)^{ C_{7} \left(
E^{\frac{193}{4}+} g^{\frac{225}{8}+}(M) \right) } \bar{r}
\end{array}
\label{Eqn:EstSizeKplus}
\end{equation}
\label{Res:BoundKplus}
\end{res}

\begin{proof}

Let $B:= \left( C_{6} E  g(M)  \right)^{ C_{7} \left(
E^{\frac{193}{4}+} g^{\frac{225}{8}+}(M) \right) } $. Assume that
(\ref{Eqn:EstSizeKplus}) fails. Let $\tilde{J}_{j_{1}}$ be the first
interval for which  $| \tilde{J}_{1} \cup .... \cup
\tilde{J}_{j_{1}} |$ exceeds $B \bar{r}$. Then $j_{1} \neq 1$,
$|\tilde{J}_{j_{1}}| \lesssim |\tilde{J}_{j_{1}-1}|
\tilde{\eta}^{-4} E^{\frac{8}{3}} g^{\frac{1}{3}}(M)$ and we have

\begin{equation}
\begin{array}{ll}
\frac{c_{1} E^{-\frac{5}{4}} \tilde{r}}{g^{\frac{5}{8}}(M)} + T_{2}
-T_{1} + (T_{2} -T_{1}) \tilde{\eta}^{-4} E^{\frac{8}{3}}
g^{\frac{1}{3}}(M) & \gtrsim |\tilde{J}_{1}| +...+ | \tilde{J}_{j_{1}} | \\
& \geq B \bar{r}
\end{array}
\label{Eqn:IneqDiffT2T1}
\end{equation}
if $[T_{1},T_{2}]:=\tilde{J}_{2} \cup ... \cup \tilde{J}_{j_{1}-1}$. Therefore by (\ref{Eqn:LowerBoundE}) and (\ref{Eqn:IneqDiffT2T1}) we have

\begin{equation}
\begin{array}{ll}
T_{2} -T_{1} & \gtrsim \frac{ \tilde{\eta}^{4} E^{-\frac{8}{3}}  B
\bar{r}}{g^{\frac{1}{3}}(M)}
\end{array}
\label{Eqn:LowerBdT1T2Diff}
\end{equation}
Moreover $ T_{1} \leq \bar{r} + \frac{c_{1} E^{-\left( \frac{5}{4} +
\right)} \bar{r}}{g^{\frac{5}{8}+}(M)} $. Therefore by
(\ref{Eqn:LowerBoundE}) we have

\begin{equation}
\begin{array}{ll}
T_{1} = O \left( \bar{r} \right)
\end{array}
\label{Eqn:UpBdT1}
\end{equation}
By (\ref{Eqn:LowerBdT1T2Diff}) and (\ref{Eqn:UpBdT1}) we have

\begin{equation}
\begin{array}{ll}
\frac{T_{2}}{T_{1}} & \geq  \left( C_{2} E^{10 +}  \left(
\frac{\tilde{\eta}}{4} \right)^{-(36+)} \right)^{4 C_{2} E^{10+}
\left( \frac{\tilde{\eta}}{4} \right)^{-(36+)}}
\end{array}
\label{Eqn:ConditionT2T1}
\end{equation}
with $C_{2}$ defined in Lemma \ref{lem:Lt4Lx12}, provided that
$C_{6}, \, C_{7} >> \max{(c_{1},C_{2})}$. Therefore we can apply
Lemma \ref{lem:Lt4Lx12} and find a subinterval $[t^{'}_{1},
t^{'}_{2} ] \subset \tilde{J}_{2} \cup....\cup \tilde{J}_{j_{1}-1}$
with $ \left| \frac{t^{'}_{2}}{t^{'}_{1}} \right| \sim E^{10+}
\tilde{\eta}^{-(36+)}$ and $\| u \|_{L_{t}^{4} L_{x}^{12}
([t^{'}_{1}, t^{'}_{2}])} \leq \frac{\tilde{\eta}}{4}$. This means
that $[ t^{'}_{1}, t^{'}_{2}] \subset [T_{1},T_{2}]$ is covered by
at most two consecutive intervals. It is convenient to introduce
$[t^{'}_{1},t^{'}_{2}]_{g}$, the geometric mean of $t^{'}_{1}$ and
$t^{'}_{2}$. We have $[t^{'}_{1},t^{'}_{2}]_{g} \sim
\tilde{\eta}^{-18} E^{5} t^{'}_{1}$. There are two cases

\begin{itemize}

\item \textit{Case 1}: $[t^{'}_{1}, t^{'}_{2}]$ is covered by one interval $\tilde{J}_{\bar{j}} = [a_{\bar{j}}, b_{\bar{j}}]$,
$2 \leq \bar{j} \leq j_{1} -1$. Then $ |\tilde{J}_{\bar{j}}| \gtrsim
\tilde{\eta}^{-(36+)} E^{10+} t^{'}_{1} $ and
$|\tilde{J}_{\bar{j}-1}| \leq t^{'}_{1}$. Therefore
$|\tilde{J}_{\bar{j}}| \gtrsim \tilde{\eta}^{-(36+)} E^{10+}
|\tilde{J}_{\bar{j}-1}|$. Contradiction with (\ref{Eqn:LowerBoundE})
and (\ref{Eqn:FractTildeJi}).

\item \textit{Case 2}: $[t^{'}_{1}, t^{'}_{2}]$ is covered by two intervals $\tilde{J}_{\bar{j}}=[a_{\bar{j}}, b_{\bar{j}}]$ and $\tilde{J}_{\bar{j}+1}=[a_{\bar{j}+1}, b_{\bar{j}+1}]$ for some $ 2 \leq \bar{j} \leq j_{1}
-2$. Then there are two subcases

\begin{itemize}

\item \textit{Case 2.a}: $b_{\bar{j}} \leq [t^{'}_{1},t^{'}_{2}]_{g} $. In this case
$|\tilde{J}_{\bar{j}+1}| \gtrsim \tilde{\eta}^{-(36+)} E^{10+}
t^{'}_{1}$ and $|\tilde{J}_{\bar{j}}| \leq \tilde{\eta}^{-(18+)}
E^{5+} t^{'}_{1}$. Therefore by (\ref{Eqn:LowerBoundE}) we have
$|\tilde{J}_{\bar{j}+1}| \gtrsim \tilde{\eta}^{-(18+)} E^{5+}
|\tilde{J}_{\bar{j}}|$. Contradiction with (\ref{Eqn:LowerBoundE})
and (\ref{Eqn:FractTildeJi}).

\item \textit{Case 2.b}: $b_{\bar{j}} \geq [t^{'}_{1},t^{'}_{2}]_{g} $. In this case by (\ref{Eqn:LowerBoundE})
$|\tilde{J}_{\bar{j}}| \gtrsim \tilde{\eta}^{-(18+)} E^{5+}
t^{'}_{1} $ and $| \tilde{J}_{\bar{j}-1} | \leq t^{'}_{1} $.
Therefore $|\tilde{J}_{\bar{j}}| \gtrsim \tilde{\eta}^{-(18+)}
E^{5+} |\tilde{J}_{\bar{j}-1}|$. Contradiction with
(\ref{Eqn:LowerBoundE}) and (\ref{Eqn:FractTildeJi}).
\end{itemize}
\end{itemize}

\end{proof}

\begin{rem}
It seems likely that we can find a better upper bound for $|K_{+}|$
than (\ref{Eqn:EstSizeKplus}) by exploiting Lemma
\ref{lem:AsympStab} in a better way. For instance we can consider
$k$ successive time intervals $\tilde{J}_{j+1}$, ....,
$\tilde{J}_{j+k}$, $k>1$ and prove an estimate like

\begin{equation}
\begin{array}{ll}
|\tilde{J}_{j+1}|+..... |\tilde{J}_{j+k}| & \lesssim |\tilde{J}_{j}|
\tilde{\eta}^{-4} E^{\frac{8}{3}} g^{\frac{1}{3}}(M)
\end{array}
\label{Eqn:StrEst}
\end{equation}
This estimate is stronger than (\ref{Eqn:FractTildeJi}). We can
probably find a smaller $B$ such that (\ref{Eqn:ConditionT2T1})
holds with $\tilde{\eta}$ substituted for something like $k
\tilde{\eta}$ and, by modifying the argument above, find a
contradiction with (\ref{Eqn:StrEst}). At the end of the process we
can probably prove global existence of smooth solutions to
(\ref{Eqn:LogQuint}) for $0< c <c_{0} $, with $c_{0}> \frac{8}{225}$
to be determined. We will not pursue these matters.

\end{rem}

\end{itemize}

\section{Proof of Lemma \ref{lem:Lt4Lx12Nrjsmall}}

Applying the Strichartz estimates and H\"older inequality

\begin{equation}
\begin{array}{ll}
\| u \|_{L_{t}^{4} L_{x}^{12} (J \times \mathbb{R}^{3})} & \lesssim E^{\frac{1}{2}} + \| u^{4} \|_{L_{t}^{1} L_{x}^{2} (J \times
\mathbb{R}^{3})} \| u g^{\frac{1}{6}}(u) \|_{L_{t}^{\infty}
L_{x}^{6} (J \times \mathbb{R}^{3}) } \| g^{\frac{5}{6}} (u) \|_{L_{t}^{\infty} L_{x}^{\infty} (J \times \mathbb{R}^{3})} \\
& \lesssim E^{\frac{1}{2}} + E^{\frac{1}{6}} g^{\frac{5}{6}}(M) \| u \|^{4}_{L_{t}^{4} L_{x}^{12} (J \times \mathbb{R}^{3})}
\end{array}
\end{equation}
Hence (\ref{Eqn:Lt4Lx12BdNrjSmall}) by (\ref{Eqn:BoundNrjM}) and a continuity argument.

\section{Proof of Lemma \ref{lem:Lt4Lx12Mass}}

Let $J^{'}=[t^{'}_{1},t^{'}_{2}] \subset J$ be such that $\| u \|_{L_{t}^{4} L_{x}^{12}(J^{'} \times \mathbb{R}^{3})}= \eta $. Then by
(\ref{Eqn:SobEmbed}) and (\ref{Eqn:UpperBdeta})

\begin{equation}
\begin{array}{ll}
\| f(u) \|_{L_{t}^{1} L_{x}^{2} (J^{'} \times \mathbb{R}^{3})} & \lesssim \| u g^{\frac{1}{6}}(u) \|_{L_{t}^{\infty} L_{x}^{6} (J^{'} \times
\mathbb{R}^{3})} \| u \|^{4}_{L_{t}^{4} L_{x}^{12} ( J^{'} \times \mathbb{R}^{3})}
\| g^{\frac{5}{6}}(u) \|_{L_{t}^{\infty} L_{x}^{\infty} (J^{'} \times \mathbb{R}^{3})} \\
& \lesssim E^{\frac{1}{6}} \eta^{4} \gfivesixM \\
& \lesssim E^{\frac{1}{2}}
\end{array}
\label{Eqn:HomEst}
\end{equation}
It is slightly unfortunate that $(2, \infty)$ is not wave admissible. Therefore we consider the admissible pair
$ \left( 2+ \epsilon, \, \frac{6(2+ \epsilon)}{\epsilon} \right)$ with $\epsilon << 1$. By the Strichartz estimates and
(\ref{Eqn:HomEst})

\begin{equation}
\begin{array}{ll}
\| u \|_{L_{t}^{2+ \epsilon} L_{x}^{\frac{6(2+ \epsilon) }{\epsilon}} (J^{'} \times \mathbb{R}^{3})  } & \lesssim \| \nabla u(t^{'}_{1})
\|_{L^{2}(\mathbb{R}^{3})} + \| u(t^{'}_{1}) \|_{L^{2} (\mathbb{R}^{3})} + \| f(u) \|_{L_{t}^{1} L_{x}^{2}
(J^{'} \times \mathbb{R}^{3})} \\
& \lesssim E^{\frac{1}{2}}
\end{array}
\label{Eqn:ControlL2Linf}
\end{equation}

Moreover let $N$ be a frequency to be chosen later. By Bernstein inequality and (\ref{Eqn:AsympF}) we have

\begin{equation}
\begin{array}{ll}
\| P_{< N} u \|_{L_{t}^{4} L_{x}^{12} (J^{'} \times \mathbb{R}^{3})} & \lesssim N^{\frac{1}{4}} |J^{'}|^{\frac{1}{4}} \| u \|_{L_{t}^{\infty} L_{x}^{6}(J^{'} \times \mathbb{R}^{3}) } \\
& \lesssim N^{\frac{1}{4}} |J^{'}|^{\frac{1}{4}} E^{\frac{1}{6}}
\end{array}
\end{equation}
Therefore

\begin{equation}
\begin{array}{ll}
\| P_{< N} u \|_{L_{t}^{4} L_{x}^{12} (J^{'} \times \mathbb{R}^{3})} & \lesssim  |J^{'}|^{\frac{1}{4}} N^{\frac{1}{4}} E^{\frac{1}{6}}
\end{array}
\label{Eqn:HighFreqBd}
\end{equation}
Let $c_{2} << 1$. Then if $ N = c_{2}^{4} \frac{\eta^{4}}{|J^{'}|
E^{\frac{2}{3}}} $ we have

\begin{equation}
\begin{array}{ll}
\| P_{\geq N} u \|_{L_{t}^{4} L_{x}^{12} (J^{'} \times \mathbb{R}^{3})} & \gtrsim \eta
\end{array}
\label{Eqn:HfLt4Lx12}
\end{equation}
and

\begin{equation}
\begin{array}{ll}
\| u \|_{L_{t}^{4} L_{x}^{12} (J^{'} \times \mathbb{R}^{3})} & \sim \| P_{\geq N} u \|_{L_{t}^{4} L_{x}^{12} (J^{'} \times \mathbb{R}^{3})}
\end{array}
\label{Eqn:EquivuHf}
\end{equation}
By (\ref{Eqn:ControlL2Linf}), (\ref{Eqn:HfLt4Lx12}) and (\ref{Eqn:EquivuHf}) we have

\begin{equation}
\begin{array}{ll}
\eta &  \sim \| P_{\geq N} u \|_{L_{t}^{4} L_{x}^{12} (J^{'} \times \mathbb{R}^{3})} \\
& \lesssim \| P_{\geq N} u \|^{\frac{2+\epsilon}{4}}_{L_{t}^{2+ \epsilon} L_{x}^{\frac{6(2+ \epsilon)}{\epsilon}}(J^{'} \times \mathbb{R}^{3})}
\| P_{\geq N} u \|^{1 - \frac{2+ \epsilon}{4}}_{L_{t}^{\infty} L_{x}^{6} (J^{'} \times \mathbb{R}^{3})} \\
& \lesssim E^{\frac{2 + \epsilon}{8}} \| P_{\geq N} u \|^{1 - \frac{2 + \epsilon}{4}}_{L_{t}^{\infty}  L_{x}^{6} (J^{'} \times \mathbb{R}^{3})}
\end{array}
\label{Eqn:Useeps1}
\end{equation}
Therefore we conclude that $ \| P_{\geq N} \|_{L_{t}^{\infty}
L_{x}^{6} (J^{'} \times \mathbb{R}^{3})} \gtrsim \eta^{2 +}
E^{-\left( \frac{1}{2} + \right)}$. Applying Proposition
\ref{prop:InvSob} we get (\ref{Eqn:Massconc}).

\section{Proof of Lemma \ref{lem:DecayPotNrjCone}}

Bahouri and Gerard ( see \cite{bahger}, p171 ) used arguments from
Grillakis \cite{grill,grill2} and Shatah-Struwe \cite{shatstruwe} to
derive an a priori estimate of the solution $u$ to the $3D$ quintic
defocusing wave equation, i.e $\partial_{tt} u - \triangle u +
u^{5}=0$. More precisely they were able to prove

\begin{equation}
\begin{array}{ll}
\int_{|x| \leq b} |u(b,x)|^{6} \, dx & \lesssim  \frac{a}{b} \left( \tilde{e}(a) + \tilde{e}^{\frac{1}{3}}(a) \right) +
\tilde{e}(b)-\tilde{e}(a) + \left( \tilde{e}(b) -\tilde{e}(a) \right)^{\frac{1}{3}}
\end{array}
\end{equation}
with

\begin{equation}
\begin{array}{ll}
\tilde{e}(t) & : = \frac{1}{2} \int_{|x| \leq t}   \left( \partial_{t} u  \right)^{2} \, dx + \frac{1}{2} \int_{|x| \leq t} | \nabla u |^{2} \,
dx + \frac{1}{6} \int_{|x| \leq t} u^{6} \, dx
\end{array}
\end{equation}
Since we apply their ideas to the potential $f$ we just sketch the proof. Given the cone $\Gamma_{+}([a,b])$ we denote by $\partial
\Gamma_{+}([a,b])$ the mantle of the cone $\Gamma_{+}([a,b])$ i.e

\begin{equation}
\begin{array}{ll}
\partial \Gamma_{+} ([a,b]) & := \left\{ (t^{'},x) \in [a,b] \times
\mathbb{R}^{3}, \, t=|x|  \right\}
\end{array}
\end{equation}
The local energy identity

\begin{equation}
\begin{array}{ll}
e(b)-e(a) & = \frac{1}{2 \sqrt{2}} \int_{\partial \Gamma_{+} ([a,b])} \left| \frac{x \partial_{t} u }{t} + \nabla u \right|^{2}  +
\frac{1}{\sqrt{2}} \int_{\partial \Gamma_{+} ([a,b])} F(u)
\end{array}
\label{Eqn:LocNrjIdent}
\end{equation}
results from the integration of the identity $\partial_{t} u \left( \partial_{tt} u - \triangle u + f(u) \right) = 0 $ on the cone
$\Gamma_{+}([a,b])$. We have \cite{shatstruwe0}

\begin{equation}
\begin{array}{l}
\partial_{t} \left( \frac{t}{2} (\partial_{t} u)^{2} + \frac{t}{2} |\nabla u |^{2} +
(x. \nabla u ) \partial_{t} u  + t F(u)  + u \partial_{t} u  \right)  \\ \\
- \dive \left( t \nabla u \partial_{t} u + (x.\nabla u) \nabla u - \frac{|\nabla u|^{2} x}{2} + \frac{ (\partial_{t} u)^{2} x}{2} - x F(u) + u
\nabla u \right) + u f(u) - 4 F(u) \\ =0
\end{array}
\end{equation}
Integrating this identity on $\Gamma_{+}([a,b])$ we have

\begin{equation}
\begin{array}{ll}
X(b) - X(a) + Y(a,b) & = \int_{\Gamma_{+}([a,b])} 4 F(u) -uf(u)
\end{array}
\end{equation}
with

\begin{equation}
\begin{array}{ll}
X(t) & := \int_{|x| \leq t} \frac{t}{2} ( \partial_{t} u )^{2} + \frac{t}{2} |\nabla u|^{2} + (x. \nabla u) \partial_{t} u + t F(u) + u
\partial_{t} u
\end{array}
\end{equation}
and

\begin{equation}
\begin{array}{ll}
Y(a,b) & := -\frac{1}{\sqrt{2}} \int_{\partial \Gamma_{+}([a,b])} \left[
\begin{array}{l}
\frac{t}{2} (\partial_{t} u)^{2} + \frac{t}{2} |\nabla u|^{2} + (x. \nabla u)
\partial_{t} u + t F(u) + u \partial_{t} u + t \frac{\nabla  u .x}{|x|} \partial_{t} u + \frac{ | x. \nabla u  |^{2}}{|x|} \\
- \frac{|\nabla u|^{2}}{2} |x| + \frac{(\partial_{t} u)^{2} |x|}{2} - |x| F(u) + u \frac{\nabla u .x}{|x|}
 \end{array}
 \right]
 \end{array}
\end{equation}
In fact (see \cite{shatstruwe}) we have

\begin{equation}
\begin{array}{ll}
X(t) & = \int_{|x| \leq t} t \left[ \frac{1}{2} (\partial_{t} u)^{2} + \frac{1}{2} \left| \nabla u + \frac{u x}{|x|^{2}} \right|^{2}   \right] +
\partial_{t} u (x.\nabla u + u) + t F(u) - \int_{|x|=t} \frac{u^{2}}{2}
\end{array}
\end{equation}
Since $t=|x|$ on $\partial \Gamma_{+} ([a,b])$ we have

\begin{equation}
\begin{array}{ll}
Y(a,b) & = -\frac{1}{\sqrt{2}} \int_{ \partial \Gamma_{+}([a,b])} |x| ( \partial_{t} u )^{2} + 2 (x.\nabla u) \partial_{t} u + u \partial_{t} u
+ \frac{(x. \nabla u)^{2}}{|x|} + u \frac{\nabla u.x}{|x|}
\end{array}
\end{equation}
and after some computations (see \cite{shatstruwe})

\begin{equation}
\begin{array}{ll}
Y(a,b) & = -\frac{1}{\sqrt{2}} \int_{\partial \Gamma_{+}([a,b])} \frac{1}{t} \left( t \partial_{t} u + (\nabla u.x) + u  \right)^{2} +
\int_{|x|=b} \frac{u^{2}}{2} - \int_{|x|=a} \frac{u^{2}}{2}
\end{array}
\end{equation}
Therefore if

\begin{equation}
\begin{array}{ll}
H(t) & :=  \int_{|x| \leq t} t \left[ \frac{1}{2} (\partial_{t} u)^{2} + \frac{1}{2} \left| \nabla u + \frac{u x}{|x|^{2}} \right|^{2}   \right]
+ \partial_{t} u (x.\nabla u + u) + t F(u)
\end{array}
\end{equation}
then

\begin{equation}
\begin{array}{ll}
H(b) - H(a) & = \frac{1}{\sqrt{2}} \int_{\partial \Gamma_{+}([a,b])} \frac{1}{t} \left( t \partial_{t} u + \nabla u.x + u \right)^{2} +
\int_{\Gamma_{+}([a,b])} 4 F(u) - uf(u)
\end{array}
\label{Eqn:VarH}
\end{equation}
We estimate $H(t)$, following \cite{bahger}. We have

\begin{equation}
\begin{array}{ll}
| \partial_{t} u (x.\nabla u + u) | & \leq \frac{t}{2} \left( (\partial_{t} u)^{2} +  \left| \nabla u  + \frac{ux}{|x|^{2}} \right|^{2} \right) \\
& \lesssim t  \left( (\partial_{t} u)^{2} + |\nabla u|^{2} + \frac{u^{2}}{|x|^{2}} \right)
\end{array}
\label{Eqn:PointwIneq}
\end{equation}
Therefore by (\ref{Eqn:PointwIneq}), H\"older inequality and (\ref{Eqn:AsympF})

\begin{equation}
\begin{array}{ll}
H(t) & \lesssim t  \left( e(t) + \int_{|x| \leq t} \frac{u^{2}}{|x|^{2}} \right)  \\
& \lesssim t \left(  e(t)+  ( \int_{|x| \leq t } u^{6} )^{\frac{1}{3}} \right) \\
& \lesssim t \left( e(t) +e^{\frac{1}{3}}(t)  \right)
\end{array}
\label{Eqn:EstH}
\end{equation}
Moreover by (\ref{Eqn:LocNrjIdent}), H\"older inequality and (\ref{Eqn:AsympF})

\begin{equation}
\begin{array}{ll}
\frac{1}{\sqrt{2}} \int_{\partial \Gamma_{+}([a,b])} \frac{1}{t} \left( t \partial_{t} u + \nabla u.x + u \right)^{2} & \lesssim \frac{b}{2
\sqrt 2} \int_{\partial \Gamma_{+}([a,b])} \left( \frac{\nabla u \cdot x}{t} + \partial_{t} u \right)^{2} + \frac{1}{2 \sqrt{2}} \int_{\partial
\Gamma_{+}([a,b])} \frac{u^{2}}{t^{2}} \\
& \lesssim b \int_{\partial \Gamma_{+}([a,b])} \left| \frac{x}{t} \partial_{t} u + \nabla u  \right|^{2}+  \frac{1}{2 \sqrt{2}} \left(
\int_{\partial \Gamma_{+}([a,b])} u^{6} \right)^{\frac{1}{3}} \\
& \lesssim b  \left( ( e(b) -e(a) ) + ( e(b)-e(a) )^{\frac{1}{3}} \right)
\end{array}
\label{Eqn:EstQuad}
\end{equation}
We get from (\ref{Eqn:AsympF})

\begin{equation}
\begin{array}{ll}
4 F(u) - u f(u) & \leq 0
\end{array}
\label{Eqn:PropFNeg}
\end{equation}
By (\ref{Eqn:VarH}), (\ref{Eqn:EstH}), (\ref{Eqn:EstQuad}) and
(\ref{Eqn:PropFNeg}) we have

\begin{equation}
\begin{array}{ll}
\int_{|x| \leq b} F(u) & \lesssim \frac{H(b)}{b} \\
& \lesssim \frac{ H(a)+ \frac{1}{\sqrt{2}} \int_{\partial \Gamma_{+}([a,b])} \frac{1}{t} \left( t \partial_{t} u + \nabla u.x + u
\right)^{2}}{b} \\
& \lesssim \frac{a}{b} \left( e(a) + e^{\frac{1}{3}}(a) \right) + e(b)-e(a) + \left( e(b) -e(a) \right)^{\frac{1}{3}}
\end{array}
\end{equation}

\section{Proof of Lemma \ref{lem:Lt4Lx12}}

The proof relies upon two results that we prove in the subsections.

\begin{res}
Let $u$ be a \textit{ classical solution} of (\ref{Eqn:LogQuint}). Assume that (\ref{Eqn:ApEstH2norm}) holds. Let $\eta$ be a positive number
such that (\ref{Eqn:UpperBdeta}) holds. If $ \| u \|_{L_{t}^{4} L_{x}^{12} (\Gamma_{+}(J))} \geq \eta$ then

\begin{equation}
\begin{array}{ll}
\| u \|_{L_{t}^{\infty} L_{x}^{6} (\Gamma_{+}(J))} & \gtrsim
\eta^{2+} E^{- \left( \frac{1}{2} + \right)}
\end{array}
\label{Eqn:Lt4Lx12PotCone}
\end{equation}
\label{res:Lt4Lx12PotCone}
\end{res}

\begin{res}
Let $u$ be a smooth solution to (\ref{Eqn:LogQuint}). Assume that (\ref{Eqn:ApEstH2norm}) holds. Let $\eta$ be a positive number such that

\begin{equation}
\begin{array}{ll}
\eta & \leq \min \left(1, E^{\frac{1}{18}} \right)
\end{array}
\label{Eqn:LongetaConst}
\end{equation}
Let $J=[t_{1}, t_{2}]$ be an interval such that $\left[ t_{1}, t_{1} \left(  E \eta^{-18} \right)^{4 E \eta^{-18}} \right] \subset J$. Then
there exists a subinterval $J^{'}=[t^{'}_{1}, t^{'}_{2}]$ such that $\left| \frac{t^{'}_{2}}{t^{'}_{1}} \right| = E \eta^{-18}$ and

\begin{equation}
\begin{array}{ll}
\| u \|_{L_{t}^{\infty} L_{x}^{6} ( \Gamma_{+}(J^{'}))} & \lesssim \eta
\end{array}
\label{Eqn:PotNrjCone}
\end{equation}
\label{res:PotNrjCone}
\end{res}
Let $C_{9}$ be the constant determined by $\gtrsim$ in (\ref{Eqn:Lt4Lx12PotCone}). Let $C_{10}$ be the constant determined by $\lesssim$ in
(\ref{Eqn:PotNrjCone}). We get from (\ref{Eqn:ConstSizeSubInt})

\begin{equation}
\begin{array}{ll}
\left[ t_{1}, \, t_{1} \left( E  \left( \frac{C_{9} \eta^{2+}
E^{-\left( \frac{1}{2} \right)+}}{2 C_{10}}  \right)^{-18}
\right)^{4 E \left( \frac{C_{9} \eta^{2+} E^{-\left( \frac{1}{2}
\right)+}}{2 C_{10}} \right)^{-18}} \right] & \subset \left[ t_{1},
\, C_{2} ( E^{10+} \eta^{-(36+)}  )^{4 C_{2} E^{10 +} \eta^{-(36+)}}
t_{1}
\right] \\
& \subset J
\end{array}
\end{equation}
if $C_{2} >> \max{(C_{9},C_{10})}$. Therefore, since $\frac{C_{9}
\eta^{2+} E^{- \left( \frac{1}{2} + \right)}}{2 C_{10}}$ satisfies
(\ref{Eqn:LongetaConst}) by (\ref{Eqn:EtaHyp}), we can use Result
\ref{res:PotNrjCone} and show that there exists a subinterval
$J^{'}=[t^{'}_{1}, \, t^{'}_{2}]$ such that $ \left|
\frac{t^{'}_{2}}{t^{'}_{1}} \right| \sim E^{10+} \eta^{-(36 + )}$
and

\begin{equation}
\begin{array}{ll}
\| u \|_{L_{t}^{\infty} L_{x}^{6} ( \Gamma_{+} (J^{'})) } & \leq
\frac{C_{9} \eta^{2+} E^{-\left( \frac{1}{2} + \right) } C_{10}}{2
C_{10}} \\
& \leq C_{9} \frac{\eta^{2+} E^{- \left( \frac{1}{2} + \right)}}{2}
\end{array}
\label{Eqn:Lt4Lx12Est}
\end{equation}
Now we claim that $\| u \|_{L_{t}^{4} L_{x}^{12} (\Gamma_{+}(J^{'}))} \leq \eta$. If not by (\ref{Eqn:EtaHyp}) and Result
\ref{res:Lt4Lx12PotCone} we have

\begin{equation}
\begin{array}{ll}
\| u \|_{L_{t}^{\infty} L_{x}^{6} ( \Gamma_{+} (J^{'}))} & \geq
C_{9} \eta^{2 +} E^{-\left( \frac{1}{2} + \right)}
\end{array}
\end{equation}
Contradiction with (\ref{Eqn:Lt4Lx12Est}).

\subsection{Proof of Result \ref{res:Lt4Lx12PotCone}}

We substitute $J^{'}$ for $\Gamma_{+}(J^{'})$ in (\ref{Eqn:HomEst}) to get

\begin{equation}
\begin{array}{ll}
\| f(u) \|_{L_{t}^{1} L_{x}^{2} (\Gamma_{+}(J^{'}) )} & \lesssim E^{\frac{1}{2}} \\
\end{array}
\label{Eqn:fuCone}
\end{equation}
By the Strichartz estimates (\ref{Eqn:StrCone}) on the truncated cone  $\Gamma_{+}(J^{'})$ we have

\begin{equation}
\begin{array}{ll}
\| u \|_{L_{t}^{2 + \epsilon} L_{x}^{\frac{6(2+ \epsilon)} {\epsilon}} (\Gamma_{+}(J^{'})) } & \lesssim E^{\frac{1}{2}}
\end{array}
\end{equation}
after following similar steps to prove (\ref{Eqn:ControlL2Linf}). Therefore

\begin{equation}
\begin{array}{ll}
\eta & = \| u \|_{L_{t}^{4} L_{x}^{12} (\Gamma_{+}(J))} \\
& \lesssim \| u \|^{\frac{2+ \epsilon}{4}}_{L_{t}^{2+ \epsilon} L_{x}^{\frac{6(2+ \epsilon)}{\epsilon}} (\Gamma_{+}(J^{'}))}
\| u \|^{1- \frac{2+ \epsilon}{4}}_{L_{t}^{\infty} L_{x}^{6} (\Gamma_{+}(J^{'}))} \\
& \lesssim E^{\frac{2 + \epsilon}{8}} \| u \|^{1- \frac{2+ \epsilon}{4}}_{L_{t}^{\infty} L_{x}^{6} (\Gamma_{+} (J^{'}))}
\end{array}
\label{Eqn:Useeps2}
\end{equation}
Therefore (\ref{Eqn:Lt4Lx12PotCone}) holds.

\subsection{Proof of Result \ref{res:PotNrjCone}}

By (\ref{Eqn:LongetaConst}) we have $ E \eta^{-18} \geq 1$. Let $n$ be the largest integer such that $2n \leq 4E \eta^{-18}$. This implies that
$n \geq E \eta^{-18}$. Let $A:= E \eta^{-18}$.  Now we consider the interval $\left[ t_{1},  A^{2n} t_{1} \right] \subset J$. We write $\left[
t_{1}, A^{2n} t_{1} \right]= \left[ t_{1}, \, A^{2} t_{1} \right] \cup ..... \cup \left[ A^{2(n-1)} t_{1}, \, A^{2n} t_{1}    \right]$. We have

\begin{equation}
\begin{array}{ll}
\sum_{i=1}^{n} e \left( A^{2i} t_{1} \right) -e \left( A^{2(i-1)} t_{1} \right) & \leq 2E
\end{array}
\end{equation}
and by the pigeonhole principle there exists $i_{0} \in [1,n]$ such that

\begin{equation}
\begin{array}{ll}
e \left( A^{2i_{0}} t_{1} \right) -e \left( A^{2(i_{0}-1)} t_{1} \right) & \lesssim \eta^{18}
\end{array}
\end{equation}
Now we choose $a:=A^{2(i_{0}-1)} t_{1} $ and $b \in [A^{2i_{0}-1}
t_{1}, \, A^{2 i_{0}} t_{1}]$. Let $t_{1}^{'}:=A^{2(i_{0}-1)}
t_{1}$, $t_{2}^{'}:=A^{2 i_{0}-1} t_{1}$ and $J^{'}:=[t_{1}^{'}, \,
t_{2}^{'}]$. We apply (\ref{Eqn:DecayPot}) and
(\ref{Eqn:LongetaConst}) to get

\begin{equation}
\begin{array}{ll}
\| u \|_{L_{t}^{\infty} L_{x}^{6} \left( \Gamma_{+} ( [t^{'}_{1}, \, t^{'}_{2}])  \right)} & \lesssim
\| F(u) \|_{L_{t}^{\infty} L_{x}^{1} \left( \Gamma_{+} ( [t^{'}_{1}, \, t^{'}_{2}])  \right)} \\
& \lesssim \left( E^{-1} \eta ^{18} (E+
E^{\frac{1}{3}}) +\eta^{18} + \eta^{6} \right)^{\frac{1}{6}} \\
& \lesssim \eta
\end{array}
\end{equation}

\section{Proof of Lemma \ref{lem:AsympStab} }

We have after computation of the derivative of $e(t)$

\begin{equation}
\begin{array}{ll}
\partial_{t} e(t) & \geq \int_{|x|=t} F(u) \, dS
\end{array}
\end{equation}
and integrating with respect of time

\begin{equation}
\begin{array}{ll}
\int_{I} \int_{|x| \leq t} \gu u^{6}(t^{'},x^{'}) \, dS \, dt^{'} & \lesssim E
\end{array}
\end{equation}
By using the space and time translation invariance

\begin{equation}
\begin{array}{ll}
\int_{J} \int_{|x^{'} - x| = |t^{'} - t|} \gu u^{6}(t^{'},x^{'}) \, dS \, dt^{'} & \lesssim E
\end{array}
\label{Eqn:ControlNrjSurf}
\end{equation}
Therefore (\ref{Eqn:ExplicitCompSin}), (\ref{Eqn:SobEmbed}), (\ref{Eqn:ControlNrjSurf}) and H\"older inequality

\begin{equation}
\begin{array}{ll}
\left| - \int_{J^{'}} \frac{\sin{(t - t^{'}) D}}{D} \gu u^{5}  \, dt^{'}  \right| \\
 = \left| \frac{1}{4 \pi |t - t^{'}|} \int_{|x^{'} -x| =|t^{'} - t|} \gfivesixu u^{5}  \gonesixu \, dS  dt^{'} \right| \\
\lesssim \int_{J^{'}} \frac{1}{|t - t^{'}|}  \left( \int_{|x^{'} - x| =|t^{'} - t|} u^{6} \gu dS \right)^{\frac{5}{6}}
\left( \int_{|x^{'} - x| =|t^{'} - t|} \gu  \, dS \right)^{\frac{1}{6}} \, dt^{'} \\
\lesssim g^{\frac{1}{6}}(M) \int_{J^{'}} \frac{1}{|t
-t^{'}|^{\frac{2}{3}}} \left( \int_{|x^{'}-x|=|t^{'}-t|} u^{6} g(u)
\, dS \right)^{\frac{5}{6}} \, dt^{'} \\
\lesssim \gonesixM  E^{\frac{5}{6}} \left( \int_{J^{'}} \frac{1}{|t - t^{'}|^{4}} \right)^{\frac{1}{6}}  \\
\lesssim \gonesixM  \frac{E^{\frac{5}{6}}}{dist^{\frac{1}{2}}(t,J^{'})}
\end{array}
\label{Eqn:ComputAsymp}
\end{equation}
Notice that

\begin{equation}
\begin{array}{ll}
u(t) & = u_{l,t_{i}}(t) - \int_{t_{i}}^{t} \frac{\sin{(t - t^{'})
D}}{ D} u^{5}(t^{'}) g(u(t^{'})) \, dt^{'}
\end{array}
\label{Eqn:DecompLinHom}
\end{equation}
for $i=1,2$. We get (\ref{Eqn:AsympStab}) from
(\ref{Eqn:ComputAsymp}) and (\ref{Eqn:DecompLinHom}).

\end{document}